\documentclass{article}[12pt,a4]
 \usepackage{latexsym,amssymb,amsmath,amsthm}
\usepackage{hyperref}
\newif\ifpdf
\ifx\pdfoutput\undefined
\pdffalse
\else
\pdfoutput=1
\pdftrue
\fi
\ifpdf
\fi
\makeatletter
\renewcommand\section{\@startsection {section}{1}{\z@}%
                                   {-3.5ex \@plus -1ex \@minus -.2ex}%
                                   {2.3ex \@plus.2ex}%
                                   {\bf\centering\normalsize}}
\renewcommand\subsection{\@startsection{subsection}{2}{\z@}%
                                     {-3.25ex\@plus -1ex \@minus -.2ex}%
                                     {1.5ex \@plus .2ex}%
                                     {\normalfont\normalsize\bfseries}}
\renewcommand\subsubsection{\@startsection{subsubsection}{3}{\z@}%
                                     {-3.25ex\@plus -1ex \@minus -.2ex}%
                                     {1.5ex \@plus .2ex}%
                                     {\normalfont\normalsize\bfseries}}
\renewcommand\paragraph{\@startsection{paragraph}{4}{\z@}%
                                    {3.25ex \@plus1ex \@minus.2ex}%
                                    {-1em}%
                                    {\normalfont\normalsize\bfseries}}
\renewcommand\subparagraph{\@startsection{subparagraph}{5}{\parindent}%
                                       {3.25ex \@plus1ex \@minus .2ex}%
                                       {-1em}%
                                      {\normalfont\normalsize\bfseries}}      
\makeatother
\newtheorem{theorem}{{ \bf Theorem}} \newtheorem{lemma}{Lemma}

\newtheorem{corollary}{Corollary} 
\newtheorem{condition}{Condition} \theoremstyle{definition}
\newtheorem{definition}{Definition} \theoremstyle{plain}
 \newtheorem{remark}{Remark}
\title{\large \bf PARAMETRIC  ESTIMATION OF DIFFUSION  PROCESSES  SAMPLED AT
  FIRST EXIT TIMES\footnote{International Journal of Pure and Applied
    Mathematics, Volume 7 No. 4 2003, 449-486}}
\author{Jaime A. Londo\~{n}o\\
{\it Departamento de Ciencias B\'asicas, Universidad EAFIT,} \\
{\it Carrera 49 \# 7 Sur-50, Medell\'{\i}n, COLOMBIA}\\
 e-mail: \texttt{jalondon@eafit.edu.co}
}\date{}

\begin{document}
\maketitle
\begin{abstract}
This paper introduces a family of recursively defined estimators 
of the parameters of a diffusion process. 
We use ideas of stochastic algorithms for the 
construction of the estimators.
Asymptotic consistency of these estimators and asymptotic normality
of an appropriate normalization are proved.  The results  are applied to
two examples from the financial literature; viz.,  Cox-Ingersoll-Ross'
model and the constant elasticity of variance (CEV) process illustrate
the use of the technique proposed herein.\newline
{\bf AMS Subj. Classification: }62M05, 62P05.\newline
\end{abstract}
\section{Introduction} \label{S:Introduction}
In this paper we introduce a family of recursively defined estimators
of the parameters of a diffusion process.  We assume that the process
is observed when it reaches some trigger values in a particular order.
An extensive literature exists on estimation for a continuous time
record of observations of diffusion processes (e.g. Bawasa an Rao \cite{iBlP80}).
Nelson \cite{dN90} studies the convergence of stochastic difference equations
to stochastic differential equations as the length of discrete time
intervals between observations goes to zero.  Banon \cite{gB78} proposes a
recursive kernel estimate of an initial density for a stationary Markov
process.

The techniques for the use of discretely-observed data are somewhat
different from those used for a continuous time record of
observations.  Maximum likelihood estimation can be applied to
discrete data, although most of the current theory requires the
discretely sample data to be stationary ergodic Markov chains.  See
 Billingsley \cite{PaBi61} or  Hall and Heyde \cite{pHcH80} for more complete references.
Lo \cite{aL88} derived a functional partial differential equation that
characterizes the likelihood function of a discretely sampled It\^o
process.  Likelihood based estimation is usually computationally quite
costly because an auxiliary partial differential equation must be
solved numerically for each hypothetical parameter value and each
observed state.  Duffie and Singleton \cite{dDkS93} and Gouri\'{e}roux
et al \cite{GMR93} suggested the use of
numerical methods to approximate moments.  He \cite{hH90} proposed the
use of binomial approximations.  Simulation approaches do not require
the Markov state vector to be fully observed.  However, it is difficult
to determine the magnitude of the approximation error, and in some
applications it might be numerically costly to ensure that the
approximation error is small.   Hansen and 
Scheinkman \cite{lHjS95} proposed moment
conditions suitable for use of generalized method of moment estimators
(see Hansen \cite{lH82}) based on properties of the infinitesimal generators
of stationary Markov processes.  A{\"{\i}}t-Sahalia \cite{yA96} proposed the use of
non-parametric techniques for the estimation of stationary
one-dimensional diffusions.  Duffie and Glynn \cite{dDpG96} introduced a family of
generalized method of moments estimators for continuous time Markov
processes observed at random time intervals.  They assume that the
arrival of the data has an intensity that varies with the underlying
Markov process or varies with an independent Markov process.  An
incomplete list of alternative estimation procedures includes
A{\"{\i}}t-Sahalia \cite{Yacine02}, Gallant, and Tauchen
\cite{Gallant_and_Tauchen96}, Stanton \cite{Stanton97},
Bandi and Phillips \cite{Bandi_and_Phillips99}, Chacko and Viceira
\cite{Chacko_and_Viceira99},
Singleton \cite{Singleton01}, Eraker \cite{Eraker97} and Jones \cite{Jones98}.

Diffusion processes play a fundamental role in stochastic optimal
control theory, stochastic thermodynamics, and financial economics.
We are particularly interested in models arising from the financial
literature.  Examples of these are exchange rate models (e.g. see
Froot and Obstfeld \cite{kFmO91} and  Krugman \cite{pK91}) and models of term structure of
interest rates (e.g. see Cox et al \cite{CIR85} and  Heath et al
\cite{HJM92}).

Our goal is to estimate the parameters of a diffusion process.  We
obtain results when the state space is one dimensional.  We assume
that the differential operator of the diffusion process
$(X_t,\mathcal{F}_t,\mathbf{P}_x)$ is given by
\begin{equation} \label{E:differential_operator}
        L= \frac{1}{2} \sigma^2(\cdot, \theta^*) \frac{d}{dx^2} +
                b(\cdot, \theta^*) \frac{d}{dx}
\end{equation}
where $b, \sigma^2$ satisfy some technical conditions sufficient for a
diffusion with this differential operator to exist, and $\theta^*
\in\mathbb{R}^s$ is a parameter to be estimated.  For any twice
continuous differentiable function $f$ on $\mathbb{R}$, it is known that
\begin{equation} \label{E:limit}
        Lf(x) = \lim_{U \downarrow \{x\}}
                \frac{\mathbf{E}_xf(X_{\tau_U})- f(x)}
                {\mathbf{E}_x(\tau_U)}
\end{equation}
where the limit is taken over open sets $U$ containing $x$ and
$\tau_U$ denotes the first exit time from the open set $U$.  For the
precise meaning of equation~(\ref{E:limit}) see Dynkin \cite{eD65}.
Therefore, it is natural to use moment conditions based on the
expressions in the numerator and the denominator of
equation~(\ref{E:limit}) to construct estimators of the parameters of
the diffusion.  It turns out that this approach suggests
parameterizations that are appropriate for identification of the
process from moment conditions of the previous type in a way that is
made precise in Section~\ref{S:Consistent}.  We use ideas from 
stochastic algorithms for the construction of the estimators.
References for the theory of stochastic algorithms are, for instance,
Benveniste et al \cite{Benveniste_et_al90}, Kushner and
Clark \cite{Kushner_and_Clark78},
 Duflo \cite{Duflo97}, Kushner and Yin \cite{hKgY97},
 Has'Minskii and Nevel'Son \cite{HaNe73}.  We prove that the sequence of estimators constructed
is asymptotically consistent and an appropriate normalization of them
is asymptotically normal.

Among the advantages of the technique that we propose are that we do
not require the diffusion process to be stationary, to have an
invariant probability measure or to satisfy some sort of ergodicity as
the techniques in prior works assume.  Another nice feature of the
estimation that we propose is the computational tractability for any
diffusion with continuous differentiable drifts and diffusion
coefficients.  In fact, we give a closed form for the functions we are
required to compute.  Finally, there exists an extensive literature
developed for the theory of stochastic algorithms, and so it is likely
that the ideas used there might be applied to this context.  A
particularly appealing characteristic of stochastic algorithms in the
econometrics of financial time series, as Benveniste et al \cite{BMP90} recall, is its
``generally recognized ability to adapt to variations in the
underlying systems''.  The latter could make it useful for the
analysis of high frequency data that seems not to be time homogeneous.

The paper is organized as follows: In Section~\ref{S:Construction} we
state some hypotheses that are used subsequently and define the
estimators that we propose. In Section~\ref{S:Consistent} we prove the
asymptotic consistency of these estimators.  In
Section~\ref{S:Normality} we prove that an appropriate normalization
of the sequence of estimators defined in Section~\ref{S:Construction}
is asymptotically normal.  In Section~\ref{S:examples} we show how the
theory we develop can be applied to some models of interest rates.
Namely, we consider Cox-Ingersoll-Ross' model and the constant
elasticity of variance (CEV) process.
\section{Construction of Estimators} \label{S:Construction} Our goal
in this paper is to estimate some parameters of a Markov process using
the values of the process that are known at some random times $
\tau_1, \tau_2, \tau_3,...$ and $\nu_2, \nu_3, ...$ that are related
to the process through the equations $~(\ref{E:stop1})$ and
$~(\ref{E:stop2})$.Genon-Catalot \cite{Genon-Catalot90},Genon-Catalot
and Laredo  \cite{Genon-Catalot_and_Laredo87}, Genon-Catalot and
Laredo \cite{Genon-Catalot_and_Laredo90},
Genon-Catalot et al \cite{Genon-Catalot02}, have constructed estimators for the parameters of a
diffusion, when only first hitting times are observed.  The use of a
space discretization rather than time discretization is well known in
the probabilistic context; it has been applied in algorithms of path
reconstruction.  See Kushner and Dupuis \cite{Kushner_and_Dupuis92},
Milstein \cite{Milstein98} and Milstein and Tretyakov \cite{Milstein_and_Tretyakov99}.

We assume that we have a parametric set of diffusions indexed by
$\mathbb{H}\subset\mathbb{R}^s$, where $\mathbb{H}$ is either a
compact set or $\mathbb{H}=\mathbb{R}^s$.  We define random variables
recursively (see equation$~(\ref{E:sequence})$ and
equation$~(\ref{E:seqtime})$) that depend on the data which we
observe.  Then we prove that, under some technical conditions, the
random variables defined in this way are asymptotically consistent for
the true value of the parameter when the parameter space is
one-dimensional.  When the parameter space is multidimensional we
obtain convergence to an invariant set of an ODE (Ordinary Differential
Equation).  From now on we assume that $\Omega=C(\left[0,\infty
\right))$ is the canonical space of continuous $S$-valued functions
where $S=\mathbb{R}$, $S=\left[0,\infty\right)$, or $S=(0,\infty)$
with the metrizable topology of uniform convergence on compact sets.
We denote by $\mathcal{F}_{\infty}$ the Borel $\sigma$-algebra of
$\Omega$.  For $0 \leq t < \infty$, $X_{t}$ is the coordinate mapping
process and,$\mathcal{F}_t$ the $\sigma$-algebra generated by
$X(\cdot)$ on $[0,t]$; namely $\mathcal{F}_t=\{X_s, 0\leq s \leq t
\}$.  We define the filtration $\mathcal{F}_{{t}^{+}}=
\bigcap_{\epsilon > 0} \sigma \, (X_{u} : \, \epsilon +t \geq u)$, for
$t \in \left[0,\infty\right)$.  For the construction of our estimators
we need the following condition:
\begin{condition} \label{C:regularity}
  $(X_t,\mathcal{F}_t,\mathbf{P}_{x}^{\theta})_{\theta \in
    \mathbb{H}}$ is a parametric set of diffusions with sample space
  $(\Omega,\mathcal{F}_{\infty})$ and differential operators
  $(L_{\theta})_{\theta \in \mathbb{H}}$.  For each $\theta \in
  \mathbb{H}\subset\mathbb{R}^s$, the part of the process
  $(X_t,\mathcal{F}_t,\mathbf{P}_x^{\theta})$ on $S$ is a recurrent
  strong Markov process.  Also
\begin{gather*}
  \sup_{x \in (a,b)} \mathbf{E}^{\theta}_{x} \, {\tau}^{(a,b)} <
  \infty \qquad (a,b) \subset S
\end{gather*}
where $\tau^{(a,b)}$ is the first exit time of the open set $(a,b)$.
We assume that the functions $\theta \mapsto \mathbf{E}^{\theta}_x Z$
and $\theta \mapsto L_{\theta}f$ are Borel measurable for
all $x \in S$, $Z$ a random variable defined on $\Omega$, and $f \in
C^2(S)$.
\end{condition}
See  Dynkin \cite{eD65,eD60} for a definition of recurrence, and of part of a
process.  In order to guarantee that a given non-negative second order
differential operator $L$ defined on $C^{2}(\mathbb{R})$ is the
differential operator of a diffusion process, it is customary to
impose:
\begin{condition} \label{C:existence}   
  $L$ is a non-negative second order differential operator with
  measurable drift coefficient $b$ and \emph{continuous} diffusion
  coefficient $\sigma^2$ that is \emph{uniformly Lipschitz continuous}
  and satisfies $\emph{either}$ of the following two properties:
 \begin{enumerate}
 \item There exists $c>0$ such that
\begin{equation*}
  \sigma^2(x) \geq c \qquad for \, x \in \mathbb{R}
\end{equation*}

  \item $\sigma^2$ is a twice continuously differentiable function such that
  the second derivative $d^2 \sigma^2/dx^2$ is bounded on
  $\mathbb{R}$.
\end{enumerate}
\end{condition}
If $L$ is as in Condition $\ref {C:existence}$ then there exists a
diffusion process ($X_t$,$\mathcal{F}_t$,$\mathbf{P}_{x}$) whose
differential operator is $L$.  (See  Kunita \cite{hK90}, Corollary 4.2.7.)

Throughout the rest of the paper we shall assume that $(X_t,
\mathcal{F}_t, \mathbf{P}_{x}^{\theta})_{\theta \in \mathbb{H}}$ is a
parametric set of one-dimensional diffusions with sample space
$(\Omega, \mathcal{F}_{\infty})$ that satisfies Condition
\ref{C:regularity} and with differential operators
$(L_{\theta})_{\theta \in \mathbb{H}}$ that satisfy Condition
\ref{C:existence}.  We assume $\theta^* \in \mathbb{H}$ is a fixed
constant and $\mu$ is a probability measure supported on $S$. We
denote by $(X_t, \mathcal{F}_t, \mathbf{P})$ the Markov process with
initial probability measure $\mu$ and probability for paths starting
at $x\in S $, $\mathbf{P}^{\theta^*}_x$; Namely $\mathbf{P}=\int
\mathbf{P}_x^{\theta^*}{d\mu(x)}$.
 
We assume that we have a finite set $\mathtt{D} =\{d_{1}, \dots
,d_{s}\} \subset S$ where $d_{1} < \dots < d_{s}$.  $\mathtt{D}$ is a
set of states where the process can be observed.  We assume that the
data arrival process is given by the following sequence of
$~(\mathcal{F}_{t^+})$ stopping times:
\begin{alignat}{2} \label{E:stop1}
  \tau^{\mathtt{D}}_1&= \inf \{t \geq 0 \mid X_{t}
  \in \mathtt{D}\} \\
  \tau^{\mathtt{D}}_{n+1}&= \inf \{ t \geq \tau_{n} \mid X_{t} \in
  \mathtt{D} \setminus \{X_{\tau_{n}}\} \}\qquad \text{for $n >1$}
  \notag
\end{alignat}
We suppress $\mathtt{D}$ in what follows.  Under the hypothesis of
recurrence we observe that $(\tau_{n})$ is a sequence of finite
stopping times.

Let $\{ U_d\}_{d \in \mathtt{D}}$ be a finite set of disjoint open
connected sets of $S$ such that $U_d \cap \mathtt{D} = \{ d \}$ for
any $d \in \mathtt{D}$. The boundary points of $U_d$ comprise a set of
states where the process can be observed given that the process has
reached the point $d \in \mathtt{D}$.  Let $\mathtt{D}_{r},
\mathtt{D}_l \colon \mathtt{D} \mapsto S \cup \{\infty, -\infty\}$ be
the functions satisfying $U_d=(\mathtt{D}_l(d),\mathtt{D}_r(d))$.  We
define $\eta^f \colon\mathbb{H}\times \mathtt{D} \mapsto \mathbb{R}$
by the formula
\begin{equation}
  \eta^f(\theta,x) = \mathbf{E}^{\theta}_x \,
  f(X_{\tau_{\mathtt{D}_r(x)} \wedge \tau_{\mathtt{D}_l(x)}}) \qquad
  \text{for } x \in \mathtt{D} \label{E:average}
\end{equation}
where $f \colon S \mapsto \mathbb{R}$ is a twice continuous
differentiable function.  Let
\begin{equation} \label{E:stop2}
  \nu_{n+1}=\inf\{t \geq \tau_{n} \mid X_t \notin U_{X_{\tau_n}}\}
  \qquad \text{for } n \geq 1
\end{equation}
We observe that $(\nu_n)_{n \geq 2}$ is a sequence of
$(\mathcal{F}_{t^+})$ stopping times.  We define $V^f$ by the formula:
\begin{equation} \label{E:deltaaverage}
  V^f(\theta,x,y)=f(y)-\eta^f (\theta,x)
\end{equation}
From now on we shall assume that the set of parameters
$\mathbb{H}=\mathbb{R}$ or $\mathbb{H}$ is the constrain set of
parameters $\mathbb{H}=\left\{\theta\colon a_i\leq\theta^i\leq
  b_i\right\}$, $-\infty <a_i<\theta^i<b_i<\infty$ for $1\leq i\leq s$
where $\theta^i$ denotes the $i$th component of $\theta$. It is
customary in the theory of stochastic algorithms to consider a
parameter set that is assumed to be  compact, due to the fact
that useful parameter values in applications are confined by
constrains of physics or economics to some constrain set.  The
constrain set mentioned above is one of such possibilities.  Other
alternatives can be consider.  See  Kushner and Yin \cite{hKgY97} or the discussion on
Section \ref{S:Consistent}.  If $H
\colon\mathbb{H}\times\mathtt{D}\to\mathbb{R}^s$ is a measurable map,
we define a sequence of estimators by the recursive relation
\begin{equation}\label{E:sequence}
  \Theta_{n+1}= \Pi_{\mathbb{H}}\left[\Theta_{n} - \gamma_nH(\Theta_{n}, X_{\tau_n})
  V^f(\Theta_{n}, X_{\tau_n}, X_{\nu_{n+1}})\right]
\end{equation}
for $n \geq 1$, where $\Theta_1$ is a bounded random variable taking
values in $\mathbb{H}$, $\Pi_{\mathbb{H}}$ is the projection onto
$\mathbb{H}$, and $(\gamma_n)$ is a decreasing sequence of positive
numbers with $\gamma_n\downarrow 0$.  In particular, $\Theta_{1}$ can
be a constant and $\Pi_{\mathbb{H}}=id$ when $\mathbb{H}=\mathbb{R}$.
The meaning of equation$~(\ref{E:sequence})$ is well known in the
theory of stochastic algorithms. A noise corrupted observation
$Y_n=H(\Theta_{n}, X_{\tau_n}) V^f(\Theta_{n}, X_{\tau_n},
X_{\nu_{n+1}})$ of a vector valued function $\bar{g}(\cdot)$ is taken,
whose root $\theta^*\in\mathbb{H}$ we are seeking.  Actually, one
observes values of the form $Y_n=g(\theta_n,X_{\tau_n})+\delta M_n$
where $\delta M_n$ has the property that $\mathbf{E}\left[\delta
  M_n\mid Y_i, \delta M_i, i <n\right]=0$.  Loosely speaking,
 $Y_n$ is an ``estimator'' of $\bar{g}(\cdot)$ in the sense that
$\bar{g}(\theta)=\lim_m (1/m)\sum_{i=1}^{m}g(\theta,X_{\tau_n})$ where
$\bar{g}(\cdot)$ is a function based on moment conditions of the type
defined by equation$~(\ref{E:average})$.  The sequence $(\gamma_n)$ is
chosen to provide an implicit average of the iterates.

In a similar way if $g \colon [0,\infty) \mapsto \mathbb{R}$ is a
measurable map, we define $\widetilde{\eta}^g \colon \mathbb{H} \times
\mathtt{D} \mapsto \mathbb{R}$, $\widetilde{V}^g \colon \mathbb{H}
\times \mathtt{D} \times [0,\infty)\mapsto \mathbb{R}$ by the
formulas:
\begin{gather}
  \widetilde{\eta}^g(\theta,x) = \mathbf{E}^{\theta}_x \,
  g(\tau_{\mathtt{D}_r(x)} \wedge \tau_{\mathtt{D}_l(x)})
        \label{E:timeaverage}\\
        \widetilde{V}^g(\theta,x,y)=g(y)-\widetilde{\eta}^g (\theta,x)
        \label{E:timedeltav}
\end{gather}
where we assume that $\{U_d\}_{d \in \mathtt{D}}$ is chosen in such a
way that $\widetilde{\eta}^g < \infty$.  For instance, under the
assumptions of Condition $\ref{C:regularity}$ it is enough to assume
that the sets $\{U_d\}_{d \in \mathtt{D}}$ have compact closure.  As
before, if
$\widetilde{H}\colon\mathbb{H}\times\mathtt{D}\to\mathbb{R}^s$ is a
measurable map, we define a sequence $(\widetilde{\Theta}_n)$ of
estimators by the recursive relation
\begin{equation} \label{E:seqtime}
 \widetilde{\Theta}_{n+1}=  \Pi_{\mathbb{H}}\left[\widetilde{\Theta}_{n} - \gamma_n
  \widetilde{H}(\widetilde{\Theta}_{n}, X_{\tau_n}) \widetilde{V}^g
  (\widetilde{\Theta}_n, X_{\tau_n}, \nu_{n+1}-\tau_n)\right]
\end{equation}
for $n \geq 1$, where $\widetilde{\Theta}_1$ is a bounded random
variable taking values in $\mathbb{H}$.  Remarks similar to the ones
done for the meaning of equation$~(\ref{E:sequence})$ hold for
equation$~(\ref{E:seqtime})$.

Instead, of using stochastic algorithms we could estimate the true
parameter of the process by trying to minimize the sum of the squares
\[
Q_n(\theta)=\sum_{k=1}^n\left[f(X_{\nu_{n+1}})-\mathbf{E}^{\theta}\left[(f(X_{\nu_{n+1}})\mid\mathcal{F}_{\tau_n}\right]\right]^2=\sum_{k=1}^n\left[f(X_{\nu_{n+1}})-\eta^f(\theta,X_{\tau_{n}})\right]^2
\]
with respect to $\theta$.  When the parameter space is unconstrained,
the estimates will be taken to be the solution of the system
\[
\frac{\partial Q_n(\theta)}{\partial\theta_i}=0, \qquad\text{for
  }i=1,\cdots,s.
\]
The ``conditional square'' approach goes back to Klimko and
Nelson \cite{Klimko_and_Nelson78} among others.  See also
 Hall and Heyde \cite{pHcH80}.  The methodological reason for which we choose to work
with stochastic algorithms instead, it is its recognized ability to
adapt to variations of the underlying system, as well as its ability
to process data sequentially as they are observed.  In the next section we will find conditions that are
sufficient for the sequence of random variables defined by
equations$~(\ref{E:sequence})$ and$~(\ref{E:seqtime})$ to be
asymptotically consistent.
\section{Consistent Estimation} \label{S:Consistent}
The next theorem is used to prove asymptotic consistency, in the case
of a one-dimensional parameterization, for the sequence of random
variables defined in equations~(\ref{E:sequence})
and~(\ref{E:seqtime}).  Compare with Theorem 7.1 from  Kushner and
Yin \cite{hKgY97}.\\
\begin{theorem}[A Robbins-Monro algorithm]\label{T:Monro}
  Let $( \Omega, \mathcal{F}, \mathbf{P})$ be a probability space, 
  $(\mathcal{F}_n)$ be a filtration of sub-$\sigma$-algebras of
  $\mathcal{F}$,  $\mathtt{D} \subset \mathbb{R}$ be a finite set, 
  and $(X_n, Y_n, \varTheta_n)_{n \in \mathbb{N}}$ be a sequence of
  real valued $(\mathcal{F}_{n})$ adapted random variables where $X_n$
  takes values in $\mathtt{D}$.  Let $\varTheta_n$ be defined by the
  following recursive relation:
\begin{equation} \label{E:Recursion}
  \varTheta_{n+1}= \varTheta_n - \gamma_n H(\varTheta_n,
  Y_n)V(\varTheta_n, Y_n,X_{n+1})
\end{equation}
where $V \colon \mathbb{R} \times \mathtt{D} \times \mathbb{R} \mapsto
\mathbb{R}$, $H \colon \mathbb{R} \times \mathtt{D} \mapsto \{1, -1\}$
are measurable functions, $(\gamma_n)$ is a decreasing sequence of
positive numbers and $\mathbf{E}(\|\varTheta_1\|^2) < \infty$ .  We
assume that the following hypotheses $A_1$, $A_2$, $H_1$, $H_2$,
and $H_3$ are satisfied:\\
$\mathrm{A_1}$ There exist a measurable function $\overline{V} \colon
\mathbb{R} \times \mathtt{D} \mapsto \mathbb{R}$ such that
\begin{equation} \label{E:A1}
  \mathbf{E}(V(\varTheta_n,Y_n,X_{n+1}) \mid \mathcal{F}_n)=
  \overline{V}(\varTheta_n,Y_n)
\end{equation}
$\mathrm{A_2}$ There exist a positive measurable function $S^2 \colon
\mathbb{R} \times \mathtt{D} \mapsto [0,\infty)$ such that
\begin{equation} \label{E:A2}
  \mathbf{E}(V^2(\varTheta_n,Y_n,X_{n+1}) \mid \mathcal{F}_n)=
  S^2(\varTheta_n,Y_n)
\end{equation}
$\mathrm{H_1}$ There exist $\theta^* \in \mathbb{R}$ such that for any
$d \in \mathtt{D}$ and $\theta \in \mathbb{R}$
\begin{equation} \label{E:H1a}
  (\theta-\theta^*) H(\theta,d)\overline{V}(\theta,d) \geq \, 0
\end{equation}
and there exists an increasing sequence of positive integers $(n_k)_{k
  \in \mathbb{N}}$ such that for any $\varepsilon >0$
\begin{equation} \label{E:H1b}
  \liminf _k \inf_{\substack{\varepsilon \, \leqslant |\theta-
      \theta^*| \,}} \mathbf{E} (        (\theta-\theta^*)
  H(\theta,Y_{n_k})\overline{V}(\theta,Y_{n_k}) )>\, 0
\end{equation}
$\mathrm{H_2}$ There exist $K>0$ such that
\begin{equation} \label{E:H2}
  S^2(\theta,d) \leq K(1 + (\theta - \theta^*)^2) \qquad \text{ for
    all } \theta \in \mathbb{R} ,\quad d \in \, \mathtt{D}
\end{equation}
$\mathrm{H_3}$ The sequence $(\gamma_n)$ of positive numbers satisfies
\begin{equation} \label{E:H3}
  \sum_{n_k} \, \gamma_{n_k} =\infty, \qquad \sum_{n} \, \gamma^2_n <
  \infty
\end{equation}

Then the sequence $(\varTheta_n)$ converges almost surely to
$\theta^*$.
\end{theorem}
 
Typically, a family $(X_t, \mathcal{F}_t,
\mathbf{P}_{x}^{(\lambda_1,\lambda_2)})_{(\lambda_1,\lambda_2) \in
  \mathbb{H}_1\times\mathbb{H}_2}$ of scalar diffusions with
differential operators
$(L_{(\lambda_1,\lambda_2)})_{(\lambda_1,\lambda_2)\in
  \mathbb{H}_1\times \mathbb{H}_2}$ is given, where
$\mathbb{H}_i\subset\mathbb{R}^{s_i}$, $i=1,2$, are compact subsets
or $\mathbb{H}_1\times \mathbb{H}_2=\mathbb{R}^{s_1}\times\mathbb{R}^{s_2}$, and
\[
L_{(\lambda_1,\lambda_2)}=\frac{1}{2}\sigma^2(x,\lambda_1)\frac{d^2}{dx^2}+
b(x,\lambda_2)\frac{d}{dx}
\]

It turns out that the sampling structures hinted by equations
(\ref{E:sequence}) and (\ref{E:seqtime}) suggest that it is more
natural to assume a parameterization defined by the indexed family of
differential operators
\begin{equation}\label{E:NewParameter}
L_{(\lambda_1^{\prime},\lambda_2^{\prime})}=\frac{1}{2}\sigma^2(x,\lambda_1^{\prime})\frac{d^2}{dx^2}+(b/\sigma^2)(x,\lambda_2^{\prime})\sigma^2(x,\lambda_1^{\prime})\frac{d}{dx}
\end{equation}
where $\lambda_i^{\prime}\in \mathbb{H}_i^{\prime}$ are compact subsets of
$\mathbb{R}^{s_i}$, $i=1,2$ or $\mathbb{H}_1^{\prime}\times
\mathbb{H}_2^{\prime}=\mathbb{R}^{s_1}\times\mathbb{R}^{s_2}$ and
where $\sigma^2$, $b/\sigma^2$ are parameterizations of the diffusion and the
ratio between the drift and the diffusion respectively.  See equations (\ref{E:calc_mom})
and (\ref{E:avwaiting2}).  It is often the case that the latter
parameterization defines an equivalent problem to the former
parameterization, at least as estimation is concerned.  Indeed, 
according to  It{\^o} and McKean \cite{Ito_and_McKean74} that borrows a phrase of W. Feller, the
expression in the numerator of
equation$~(\ref{E:differential_operator})$ defines a road map, i.e. it
tells what routes the particle is permitted to travel, and the
expression at the bottom of equation$~(\ref{E:differential_operator})$
defines the speed of the diffusion. Using Feller's terminology,
$\lambda_2^{\prime}$ identifies the ``road map'', and $\lambda_1^{\prime}$
identifies the ``speed'' of the diffusion when the ``road map'' is known.  In this paper
we should adopt the latter approach.  Corollary \ref{C:scalecons}
below is used to estimate the parameter that identifies the ratio
between the drift and the diffusion when the parameter space used to
identify this ratio is $\mathbb{R}$.  Theorem \ref{T:consistencyMulty}
is used when a compact subset of $\mathbb{R}^{s_1}$ is used as the
parameter space that identifies this ratio.  Let us observe that neither case  requires the
parameter(s) that identifies the diffusion to be known.  Thus, it is
possible to assume that the latter parameter is known, when the
estimation of the former parameters of  the diffusion is made.

For the next two corollaries, let us assume that the parameter space
is one dimensional.
\begin{corollary} \label{C:scalecons}
  Let $\sigma^2 \colon S \times \mathbb{R} \mapsto [0,\infty)$, and $b
  \colon S \times \mathbb{R} \mapsto \mathbb{R}$ be defined by the
  formula
\begin{equation}
  L_{\lambda}= \frac{1}{2} \sigma^2(\cdot,\lambda) \frac{d^2}{dx^2} +
  b(\cdot,\lambda) \frac{d}{dx} \qquad \text{ for } \lambda \in
  \mathbb{R}
\end{equation}
where $b(\cdot,\lambda)$, $\sigma^2(\cdot,\lambda)$,
$(b/\sigma^2)(\cdot,\lambda) \in C(\overline{S}) \cap C^2(S)$ for any
$\lambda \in \mathbb{R}$.  We assume that $\partial/\partial \lambda
(b/ \sigma^2)(x,\lambda)$ exists and is nowhere zero for $(x,\lambda)
\in
\cup_{d \in \mathtt{D}} U_d \times \mathbb{R}$. \\
Let $H \colon \mathbb{R} \times \mathtt{D} \to \{-1, 1\}$ be defined
by the formula
\begin{equation} 
  H(\lambda,d)= \mathbf{1}_{(0,\infty)}(\partial /\partial
  \lambda(b/\sigma^2)(d,\lambda))-
  \mathbf{1}_{(-\infty,0)}(\partial/\partial \lambda
  (b/\sigma^2)(d,\lambda))
\end{equation}

If $\lambda^* \in \mathbb{R}$ is a fixed number, then the sequence of
random variables defined by equation~(\ref{E:sequence}), converges
almost surely $\mathbf{P}^{\lambda^*}_x$ to $\lambda^*$ for any $x \in
S$, where $\eta^f$ and $V^f$ are defined as in
equations~(\ref{E:average}) and ~(\ref{E:deltaaverage}), and $f=id$ is
the identity on $\mathbb{R}$.
 \end{corollary}

A few words are needed to review the hypotheses from Corollary
\ref{C:scalecons}.  If the drift is zero the sampling scheme defined
by equation~(\ref{E:sequence}) can not be used.  In fact, only data
obtained using the sampling scheme defined by
equation~(\ref{E:timedeltav}) would provide any information.  See
Corollary \ref{C:timecons} below.  Within the framework proposed this
is indeed natural. If the drift is zero, it is conceivable that only
the times between hits of the grids and the end points of the
surrounding intervals should provide any information.  If
$b/\sigma^2(d,\cdot)$, $d\in\mathtt{D}$ are strictly monotone
functions around an interval containing the ``true'' parameter, then
it is possible to define a new parameterization that complies with the
hypothesis of the previous corollary and allows us to identify the
parameter at least from a small interval.  Also, Theorem
\ref{T:consistencyMulty} can be used whenever $b/\sigma^2(d,\cdot)$,
$d\in\mathtt{D}$ are not strictly monotone.

Corollary \ref{C:timecons} below is used to estimate the parameter
that identifies the diffusion, when the parameter space to identify
this diffusion is $\mathbb{R}$.  See Theorem \ref{T:consistencyMulty}
for estimation of parameters used to identify the diffusion term for a
multidimensional setting for the parameter space.  It
is assumed that the vector of parameter(s) that identifies the ratio
between the drift and the diffusion is known.  The previous assumption
can be made in lieu of Corollary \ref{C:scalecons} or Theorem
\ref{T:consistencyMulty} in conjunction with the remarks made right
after Theorem \ref{T:Monro}.

\begin{corollary} \label{C:timecons}
  Let $\sigma^2 \colon S \times \mathbb{R} \mapsto [0,\infty)$, and $b
  \colon S \times \mathbb{R} \mapsto \mathbb{R}$ be defined by the
  formula
\begin{equation}
  L_{\varsigma}= \frac{1}{2} \sigma^2(\cdot,\varsigma)
  \frac{d^2}{dx^2} + b(\cdot,\varsigma) \frac{d}{dx} \qquad \text{ for
    } \varsigma \in \mathbb{R}
\end{equation}
where $b(\cdot,\varsigma)$, $\sigma^2(\cdot,\varsigma)$,
$b/\sigma^2(\cdot,\varsigma) \in C(\overline{S}) \cap C^2(S)$ for any
$\varsigma \in \mathbb{R}$.  We assume:
\begin{enumerate}
\item There exists a function $s \colon S \mapsto \mathbb{R}$, that
  does not depend on $\varsigma$, such that
\begin{equation}\label{E:independence}
  \frac{b(x,\varsigma)}{\sigma^2(x,\varsigma)}=s(x) \qquad \text{for
    any } x \in S, \, \varsigma \in \mathbb{R}
\end{equation}\label{P:scale}
\item There exist $\sigma_0 \colon \mathbb{R} \mapsto \mathbb{R}^+$,
  $h \colon S \mapsto \mathbb{R}$ such that
\begin{equation}\label{E:decomposition}
  \sigma^2(x,\varsigma)=\sigma_0(\varsigma)h(x) \qquad \text{ for any
    } x \in S, \, \varsigma \in \mathbb{R}
\end{equation} 
where we assume that $\sigma_0$ is a strictly increasing function that
is differentiable and
\begin{equation} \label{E:rate}
  \lim_{s \to \infty} \inf_{ |\varsigma| \geq s} |\varsigma
  |\sigma_0(\varsigma) >0
\end{equation}
\label{P:separability}
\end{enumerate}
If $\varsigma^* \in \mathbb{R}$ is a fixed number, then the sequence
of random variables defined by equation~(\ref{E:seqtime}) converges
almost surely $\mathbf{P}^{\varsigma^*}_x$ to $\varsigma^*$ for any $x
\in S$, where $\widetilde{\eta}^g$ and $\widetilde{V}^g$ are defined
as in equations~(\ref{E:timeaverage}) and ~(\ref{E:timedeltav}) for
$g$ the identity on $\mathbb{R}^+$ and $\widetilde{H} \colon
\mathbb{R} \times \mathtt{D} \mapsto \{ -1, 1\}$ is the constant
function equal to $1$.
\end{corollary}
Let us review the hypotheses of Corollary \ref{C:timecons}.
Equation~(\ref{E:independence}) is natural under the assumptions made
on the parameterization.  See the remarks made after Theorem
\ref{T:Monro}.  The factorization of equation~(\ref{E:decomposition})
often arises in applications. The latter assumption is used to prove
monotonicity of the function defined by equation~(\ref{E:timecons1}).
The assumption made on equation~(\ref{E:rate}) can be made without any
loss of generality.

In order to illustrate the use of stochastic algorithms for the
problem of estimation in the multidimensional case (for the parameter
space), we make use of the standard theorem of convergence for
truncated stochastic algorithms with correlated noise with step size
going to zero.
Assume a constrain multidimensional parameter space
$\mathbb{H}=\{\theta\colon a_i\leq\theta^i\leq b_i\}$,
$-\infty< a_i<\theta^i<b_i<\infty$ for $1\leq i\leq s$.  For $\theta\in \mathbb{H}$, define the set $C(\theta)$ as
follows.  For $\theta\in\mathbb{H}^{0}$, the interior of $\mathbb{H}$, $C(\theta)$
contains only the zero element; for $\theta\in\partial\mathbb{H}$, the
boundary of $\mathbb{H}$, let $C(\theta)$ be the infinite convex cone generated
by the outer normals at $\theta$ on the faces on which $\theta$ lies.
Given a continuous $g\colon\mathbb{R}^s\mapsto\mathbb{R}^s$ the
\emph{projected} ODE of $\dot{\theta}=g(\theta)$ is defined to be
\[
\dot{\theta}=g(\theta)+z,\qquad \theta(t)\in -C(\theta(t))
\]
where $z(\cdot)$ is the projection or constrain term, the minimum term
needed to keep $\theta(\cdot)$ in $\mathbb{H}$.

\begin{theorem}\label{T:consistencyMulty}
  Let $( Y_n, \varTheta_n)_{n \in \mathbb{N}}$ be a sequence of
  $(\mathcal{F}_{\tau_n})$ adapted measurable maps where $Y_n \colon
  (\Omega, \mathcal{F}_{\tau_n}) \mapsto
  (\mathbb{R},{\mathcal{B}}(\mathbb{R}))$, and $\varTheta_n \colon
  (\Omega,\mathcal{F}_{\tau_n}) \mapsto (\mathbb{R}^s,
  \mathcal{B}(\mathbb{R}^s))$. Let $K$ a non singular $s\times s$
  matrix. Assume that $(Y_n,\varTheta_n)$ satisfies the following
  recursive relation:
\begin{equation} \label{E:MultiRecursion}
    \varTheta_{n+1}=\Pi_{\mathbb{H}}\left[\varTheta_n - \gamma_n K\triangledown\overline{V} (\varTheta_n,
    X_{\tau_n})V(\varTheta_n, X_{\tau_n} ,Y_{n+1})\right]
\end{equation}
where $\Pi_{\mathbb{H}}$ is the projection onto $\mathbb{H}$, $V \colon\mathbb{H}\times
\mathtt{D} \times \mathbb{R} \mapsto \mathbb{R}$ is a measurable
function, $\overline{V} \colon\mathbb{H}\times \mathtt{D} \mapsto \mathbb{R}$
is a twice continuous differentiable function with differential
$\triangledown\overline{V}(\cdot,d)$ for $d\in\mathtt{D}$, $K$ is an
invertible matrix and $(\gamma_n)$ is a decreasing sequence of
positive numbers.  We assume that
\begin{equation} \label{E: MA1a}
    \mathbf{E}(V(\varTheta_n,X_{\tau_n},Y_{n+1}) \mid \mathcal{F}_{\tau_n})=
        \overline{V}(\varTheta_n,X_{\tau_n})\qquad\mbox{ for }n\geq 1
\end{equation}
Moreover, it is assumed that the sequence $(\gamma_n)$,
$\gamma_n\downarrow 0$ of positive numbers satisfies
\begin{equation}\label{E:StepCondition}
\sum_{n}\gamma_{n}=\infty , \qquad\sum_{n}\gamma^2_n  <\infty .
\end{equation}
Then the sequence $(\varTheta_n)$ converges almost surely $\mathbf{P}$
to an invariant set of the projected ODE
\begin{equation}\label{E:PODE}
\dot{\theta}=-K\overline{g}(\theta)+z,\qquad \theta(t)\in
-C(\theta(t)))
\end{equation}
for
\[
\overline{g}(\theta)=\frac{1}{2}\sum_{d\in\mathtt{D}}{p_d\triangledown}\overline{V}^2(\theta,d)
\]
where $p=(p_i)$ is the left-fixed probability row vector for the
Markov chain $(X_{\tau_n},\mathcal{F}_{\tau_n})$. Indeed,
$(\varTheta_n)$ converges almost surely to a unique compact and
connected component of the set of stationary points of the
equation~(\ref{E:PODE}). If $\theta^*$ is an asymptotically stable
point of equation~(\ref{E:PODE}) and $(\varTheta_n)$ is in some
compact set in the domain of attraction of $\theta^*$ infinitely often
with probability $\geq\rho$, then $\varTheta_n\to\theta^*$ with at
least probability $\rho$.
\end{theorem}
The proof of the above theorem is a straightforward
consequence of the Theorem of convergence with probability one for the
correlated noise case for stochastic algorithms.  See for example
 Kushner and Yin \cite{hKgY97}, Theorem 6.1.1. The details of the proof are left to the
reader.

More general constrain sets can be consider. For instance, let
$q_i(\cdot), \, i=1,\cdots,p$ be continuously differentiable
real-valued functions on $\mathbb{R}^s$, with gradients $\triangledown
q_i(\cdot)$, where it is assumed that $\triangledown q_i(x)\neq 0$ if
$q_i(x)=0$ and that $\mathbb{H}=\{x\mid q_i(x) \leq 0, i=1,\cdots,p\}$
is a nonempty, compact connected set.  Define $C(x)$ to be the convex
cone generated by the set of outward normals $\left\{\triangledown
  q_i(x)\mid q_i(x)=0\right\}$.  Suppose that for each $x$ the set
$\left\{\triangledown q_i(x)\mid q_i(x)=0\right\}$ is either empty or
a linear independently set.  Then the Theorem \ref{T:consistencyMulty}
remains true with the obvious changes.  See  Kushner and Yin
\cite{hKgY97}.  Similarly,
if $\mathbb{H}$ is a $\mathbb{R}^{s-1}$ dimensional connected compact surface
with a continuous differentiable outer normal, and we define $C(x)$,
$x\in \mathbb{H}$, to be the linear span of the outer normal at $x$ then
Theorem \ref{T:consistencyMulty} still holds.  See also
 Kushner and Yin \cite{hKgY97}.  It is worth noting that the former constrain set, as
well as the mentioned in the Theorem \ref{T:consistencyMulty}, can give
rise to new stationary points of the ODE $~(\ref{E:PODE})$, but this
is the only type of singular point that can be introduced by the
constrains.  In many applications when the truncation bounds are large
enough, there is only one stationary point $\theta^*$ of the ODE
$~(\ref{E:PODE})$ that is globally asymptotically stable.  Typically,
for the kind of application we are heading,
$\bar{V}(\theta,d)=\eta(\theta^*,d)-\eta(\theta,d)$ where $\eta\colon
\mathbb{H}\times\mathtt{D}\mapsto\mathbb{R}^s$ is a twice continuous
differentiable function (on the parameter variable) and $\theta^*\in
\mathbb{H}^0$.  If
\[
\left[\begin{array}{ccc}
    \frac{\partial\eta}{\partial\theta^1}(\theta^*,d_1 )&\cdots  &\frac{\partial\eta}{\partial\theta^s}(\theta^*,d_1)  \\
    \vdots &\ddots  &\vdots  \\
    \frac{\partial\eta}{\partial\theta^1}(\theta^*,d_r)&\cdots &
    \frac{\partial\eta}{\partial\theta^s}(\theta^*,d_r)
\end{array}\right]
\]
defines an injection, then $\theta^*$ is the unique stationary point
of equation $~(\ref{E:PODE})$ in the interior of $\mathbb{H}$, at least for a
sufficiently small neighborhood of $\theta^*$.

As a application of Theorem \ref{T:consistencyMulty}, let us assume
a family of scalar diffusions $(X_t, \mathcal{F}_t,
\mathbf{P}_{x}^{(\lambda_1,\lambda_2)})_{(\lambda_1,\lambda_2) \in
  \mathbb{H}_1\times \mathbb{H}_2}$ with differential operators
$(L_{(\lambda_1,\lambda_2)})_{(\lambda_1,\lambda_2)\in
  \mathbb{H}_1\times \mathbb{H}_2}$ given by equation
$~(\ref{E:NewParameter})$ where $\mathbb{H}_i\subset\mathbb{R}^{s_i}$,
$i=1,2$ are constrain sets as the ones discussed above, for $i= 1,
2$.  It is assumed that there exist a parameter
$(\lambda^*_1,\lambda^*_2)\in \mathbb{H}_1\times \mathbb{H}_2$ such
that $(\mathbf{P}^{(\lambda^*_1,\lambda^*_2)}_x)=(\mathbf{P}_x)$.  We
consider the family of diffusions $(X_t,
\mathcal{F}_t,\mathbf{P}_{x}^{(\lambda,\lambda_2^*)})_{\lambda \in
  \mathbb{H}_1}$. Let $(\Theta_n)$ be defined by
equation$~(\ref{E:sequence})$, where the projection in taken over the
set $\mathbb{H}_1$,
$\eta^f(\lambda,x)=\mathbf{E}_x^{\lambda,\lambda^*_2} \,
f(X_{\tau_{\mathtt{D}_r(x)} \wedge \tau_{\mathtt{D}_l(x)}})$,
$V^f(\lambda,x,y)=f(y)-\eta^f(\lambda,x)$,
$\overline{V}^f(\lambda,x)=\eta^f(\lambda_1^*,x)-\eta^f(\lambda,x)$,
$K$ is a non-singular matrix and $(\gamma_n)$ is a sequence as in
equation$~(\ref{E:StepCondition})$. It follows that Theorem
$\ref{T:consistencyMulty}$ applies, and it identifies $\lambda_1^*$ if
this is the unique stationary point of the projected ODE
$~(\ref{E:PODE})$.  We observe that the computation made to obtain the
sequence $(\Theta_n)$ does not depend on the value $\lambda_2^*$. See
Appendix \ref{S:Solution1} for the computation of the algorithms.
Next, we assume that the parameter $\lambda_1^*$ is known.  (The
latter can be assumed by the previous remark.)  We consider the
collection of diffusions
$(X_t,\mathcal{F}_t,\mathbf{P}_x^{(\lambda_1^*,\lambda^{\prime})})_{\lambda^{\prime}\in\mathbb{H}_2}$.
Let $(\tilde{\Theta}_n)$ be the sequence of estimators defined by
equation$~(\ref{E:seqtime})$ where the projection is taken oven the
set $\mathbb{H}_2$,
$\eta^g(\lambda^{\prime},x)=\mathbf{E}^{\lambda_1^*,\lambda^{\prime}}\,g(\tau_{\mathtt{D}_r(x)}
\wedge \tau_{\mathtt{D}_l(x)})$, $
\tilde{V}^g(\theta,x,y)=g(y)-\tilde{\eta}^g (\theta,x)$,
$\overline{V}^g(\lambda^{\prime},x)=\tilde{\eta}^{g}(\lambda^*_2,x)-\tilde{\eta}^{g}(\lambda,x)$,
$K$ a non-singular matrix (not necessarily identical to the one used to
compute $(\Theta_n)$), then Theorem $\ref{T:consistencyMulty}$
applies, and it identifies $\lambda_2^*$ if this is the unique
stationary point of the projected ODE $~(\ref{E:PODE})$.

It is worth noting that even if only data associated with the sampling
scheme related with equation$~( \ref{E:sequence})$ is available then
at least, identification of the the ratio between the drift and the
diffusion can be made.  Also, when the dimension of the parameter
space that identifies either the diffusion or the ratio between the
drift and the diffusion are one-dimensional, Corollaries
\ref{C:scalecons} and \ref{C:timecons} can be called for the
estimation with the advantage that complete identification of the
parameter is easier.    
\section{Asymptotic Normality} \label{S:Normality}
In this section we propose a version of the central limit theorem
for the class of estimators of Theorem~\ref{T:Monro}.
  
For any stopping
time $\tau$ we denote as $\theta_{\tau}$ the measurable map
defined as $\theta_{\tau}(\omega)(\cdot)= \theta(\tau(\omega) +
\cdot)$.  We observe that
$\theta_{\tau_{n}}=\theta^{n-1}_{\tau_2}$ and $X_{\tau_{n}}=
X_{\tau_1} \circ \theta_{\tau_n}= X_{\tau_1} \circ
\theta^{n-1}_{\tau_2}$ for $n \geq 2$.

For the following theorem we assume that
$(X_t$,$\mathcal{F}_t$,$\mathbf{P}_{x}^{\theta})_{\theta \in
  \mathbb{R}}$ is a parametric family of recurrent strong Markov
processes . 
\begin{theorem} \label{T:Normality_One_dimension}
  Let $Y\colon (\Omega, \mathcal{F}_{\tau_2}) \to (\mathbb{R},
  \mathcal{B}(\mathbb{R}))$ be a measurable map that is bounded below.
  Moreover, assume that $Y \in \bigcap_{d \in \mathtt{D}}
  L^2(\mathbf{P}_d)$.  Let $Y_n$ be defined as $Y_n=Y \circ
  \theta_{\tau_{n-1}} =Y \circ \theta^{n-2}_{\tau_2}$, for $n \geq 3$
  and $Y_2=Y$.  Let $\eta \colon \mathbb{R} \times \mathtt{D} \to
  \mathbb{R}$ be the function defined as
  $\eta(\theta,d)=\mathbf{E}_d^{\theta} (Y)$.  In addition assume that
  Hypotheses $N_1$ and $N_2$
  are satisfied:\\
  $\mathrm{N_1}$ For any $d \in \mathtt{D}$ $\eta(\cdot,d)$ is a
  strictly monotone, twice continuous differentiable function with
  non-vanishing
  derivative.\\
  $\mathrm{N_2}$ There exist $L, L' > 0$ such that for any $\theta \in
  \mathbb{R}$, $d \in \mathtt{D}$
        \begin{equation} \label{E:Liptchitz}
        \mid  \eta(\theta,d) -\eta(\theta^*,d) \mid
        \leq L \mid \theta- \theta^* \mid  + L'
        \end{equation}
        Define $V \colon \mathbb{R} \times \mathtt{D} \times
        \mathbb{R} \to \mathbb{R}$ by the formula
        $V(\theta,d,y)=y-\eta(\theta,d)$.  Assume that $\varTheta^N_n$ is a
        $(\mathcal{F}_{\tau_n})$ adapted sequence of random variables
        that satisfies the recursive relation:
\begin{equation} \label{E:Normal_Recursion}
        \Theta^N_{n+1}= \Theta^N_n -\frac{1}{n}
        \frac{V(\Theta^N_n, X_{\tau_n},Y_{n+1})}{\alpha(X_{\tau_n})}
\end{equation}
where $\alpha(d)=-(\partial\eta/\partial \theta)(\theta^*,d)$ for
$d \in \mathtt{D}$ and $\mathbf{E}((\Theta_1)^2) < \infty$.\\ Then
$n^{1/2} (\Theta^N_{n}-\theta^*)$ is asymptotically normally
distributed with mean zero and variance $\sigma^2=\sum_{d \in
  \mathtt{D}}p_d Var_d(Y)/\alpha^2(d) =\sum_{d \in \mathtt{D}}p_d
\mathbf{E}_d(Y - \eta(\theta^*,d))^2/ \alpha^2(d)$ where $p=(p_i)$
is the left-fixed probability row vector of the Markov chain
$(X_{\tau_n}, \mathcal{F}_{\tau_n})$ as in
Lemma~\ref{L:Ergodicity}.
\end{theorem}


\begin{corollary} \label{C:scalenormal}
  Let $\mu$ be a probability measure $\mathbb{R}$ supported on $S$,
  $\lambda^* \in \mathbb{R}$ be a fixed constant, and let $(\Omega,
  \mathbf{P}, \mathcal{F}_{\infty})$ be the probability space where
  $\mathbf{P}=\int \mathbf{P}^{\lambda^*}_x d \mu(x)$.  Let $b,
  \sigma^2, b/\sigma^2$ and $(L_{\lambda})$ be as in Corollary~\ref{C:scalecons}.
Let $\eta \colon \mathbb{R} \times \mathtt{D}
  \to \mathbb{R}$ be defined as
  $\eta(\lambda,d)=\mathbf{E}_d^{\lambda} (X_{\nu_2})$.  We define $V
  \colon \mathbb{R} \times \mathtt{D} \times \mathbb{R} \to
  \mathbb{R}$ by the formula $V(\lambda,d,y)=y-\eta(\lambda,d)$.
  Assume that $(\Theta^N_n)$ is a $(\mathcal{F}_{\tau_n})$ adapted
  sequence of random variables which satisfies the recursive relation
\begin{equation} \label{E:Normal_Recursion_scale}
        \Theta^N_{n+1}= \Theta^N_n -\frac{1}{n}
        \frac{V(\Theta^N_n, X_{\tau_n},X_{\nu_{n+1}})}
        {\alpha(X_{\tau_n})}
\end{equation}
where $\alpha(d)=-(\partial\eta/\partial \lambda)(\lambda^*,d)$
for $d \in \mathtt{D}$ and $\Theta^N_1$ is a bounded random
variable.  Then $n^{1/2}(\Theta^N_n - \lambda^*)$ is
asymptotically normally distributed with mean zero and variance
$\sigma^2 =\sum_{d \in
  \mathtt{D}}p_d \mathbf{E}^{\lambda^*}_d (X_{\nu_2} -
\eta(\lambda^*,d))^2/\alpha^2(d)$ where $p=(p_i)$ is the
left-fixed probability row vector of the Markov chain
$(X_{\tau_n}, \mathcal{F}_{\tau_n})$ as in
Lemma~\ref{L:Ergodicity}.
\end{corollary}
\begin{corollary} \label{C:timecons_normality}
  Let $\mu$ be a probability measure on $\mathbb{R}$ supported on $S$,
  $\varsigma^* \in \mathbb{R}$ be a fixed constant, and let $(\Omega,
  \mathbf{P}, \mathcal{F}_{\infty})$ be the probability space where
  $\mathbf{P}=\int \mathbf{P}^{\varsigma^*}_x d \mu(x)$.  We assume
  that $b$, $\sigma^2$, $(L_{\varsigma})$, $s$, $\sigma_0$, and $h$
  satisfy the hypothesis of Corollary~\ref{C:timecons}.  Let
  $\tilde{\eta} \colon \mathbb{R} \times \mathtt{D} \to \mathbb{R}$ be
  defined as $\tilde{\eta}(\varsigma,d)=\mathbf{E}_d^{\varsigma}
  (\nu_2)$.  We define $\tilde{V} \colon \mathbb{R} \times \mathtt{D}
  \times \mathbb{R} \to \mathbb{R}$ by the formula
  $\tilde{V}(\varsigma,d,y)= y-\tilde{\eta}(\varsigma,d)$.  Assume
  that $(\tilde{\Theta}^N_n)$ is a $(\mathcal{F}_n)$ adapted sequence
  of random variables that satisfies the recursive relation
\begin{equation} \label{E:Normal_Recursion_timecons}
        \tilde{\Theta}^N_{n+1}= \tilde{\Theta}^N_n -\frac{1}{n}
        \frac{\tilde{V}(\tilde{\Theta}^N_n, X_{\tau_n},\nu_{n+1}-\tau_n)}
        {\tilde{\alpha}(X_{\tau_n})}
\end{equation}
where $\tilde{\alpha}(d)=-(\partial\tilde{\eta}/\partial
\varsigma) (\varsigma^*,d)$ for $d \in \mathtt{D}$ and
$\tilde{\Theta}^N_1$ is a bounded random variable.  Then
$n^{1/2}(\tilde{\Theta}^N_n - \varsigma^*)$ is asymptotically
normally distributed with mean zero and variance $\sigma^2
=\sum_{d \in \mathtt{D}}p_d \mathbf{E}^{\varsigma^*}_d (\nu_2 -
\tilde{\eta}(\varsigma^*,d))^2/\tilde{\alpha}^2(d)$ where
$p=(p_i)$ is the left-fixed probability row vector of the Markov
chain $(X_{\tau_n}, \mathcal{F}_{\tau_n})$ as in
Lemma~\ref{L:Ergodicity}.
\end{corollary}
In order to illustrate the use of stochastic algorithms, for the
problem of asymptotic normality in the multidimensional case (for the
parameter space), we use a standard theorem for the rate of
convergence for stochastic algorithms with correlated noise and
decreasing step size.  We assume a constrain multidimensional
parameter space $\mathbb{H}=\left\{\theta\colon a_i\leq\theta^i\leq
  b_i\right\}$, $-\infty <a_i<\theta^i<b_i<\infty$ for $1\leq i\leq
s$.  We assume that
$(X_t$,$\mathcal{F}_t$,$\mathbf{P}_{x}^{\theta})_{\theta \in \mathbb{H}}$ is a
parametric family of recurrent strong Markov processes with sample
space $(\Omega, \mathcal{F}_{\infty})$. Let $\theta^* \in \mathbb{H}$ be an
interior point, and $\mu$ is a probability measure supported on $S$.
We denote as $(X_t, \mathcal{F}_t, \mathbf{P})$ the Markov process
with parameter $\theta^*$ and initial probability measure $\mu$.  Let
$Y$, $Y_n$ $n \geq 2$ be defined as in Theorem
\ref{T:Normality_One_dimension}, and define $\eta\colon
\mathbb{H}\times\mathtt{D}\to\mathbb{R}$ as
$\eta(\theta,d)=\mathbf{E}_d^{\theta}(Y)$.  We assume a recursive
sequence of estimators $(\varTheta_n)$ defined as,
\begin{equation} \label{E:NormalMultiRecursion}
    \varTheta_{n+1}=\Pi_{\mathbb{H}}\left[\varTheta_n +\frac{1}{n}K\triangledown\eta(\varTheta_n,X_{\tau_n}) (Y_{n+1}-\eta(\varTheta_n,X_{\tau_n}))\right]
\end{equation}

Let $D\left(-\infty,\infty\right)$ (resp., $D\left[0,\infty\right)$)
denote the space of real-valued functions on the interval
$\left(-\infty,\infty\right)$ (resp., on $\left[0,\infty\right)$) that
are right continuous and have left hand limits, endowed with the Skorohod
topology, and $D^s\left(-\infty,\infty\right)$ its $s$-fold
product. Full descriptions and treatments of the Skorohod topology are
given in Billingsley \cite{Billingsley99} and  Ethier and Kurtz \cite{sEtK86}. Define
$U_n=\sqrt{n}(\varTheta_{n}-\theta^*)$, and let $U^n(\cdot)$ denote
the piecewise constant right continuous interpolation (with
interpolation intervals $\{1/n\}$) of the sequence $\{U_i,i\geq n\}$
on $[0,\infty)$.  Namely, if we define $t_0=0$ and
$t_n=\sum_{i=1}^{n}1/i$, we make
\begin{equation*}
  U^n(t)=U_m\text{ for }t\in\left[t_{m-1}-t_{n-1},t_m-t_{n-1}\right),\text{ and }m\geq n\geq 1
\end{equation*}
For $t\geq 0$, let $m(t)$ denote the unique value of $n$ such that
$t\in\left[t_{n-1},t_{n}\right)$, and for $t<0$ set $m(t)=1$. Define
the continuous time interpolation $W^n(\cdot)$ on
$\left(-\infty,\infty\right)$, for $n\geq 1$, by
\begin{equation}\label{E:NormalError}
W^n(t)  = \left\{ \begin{array}{c}
 \sum\limits_{i = n}^{m(t_n  + t) - 1} {\frac{1}{{\sqrt i }}\left(K\triangledown\eta(\theta^*,X_{\tau_i})(Y_{i+1}-\eta(\theta^*,X_{\tau_i})) \right)\qquad\text{for t }\geq 0}  \\ 
  \\ 
 -\sum\limits_{m(t_n  + t)}^n {\frac{1}{{\sqrt i }}\left(K\triangledown\eta(\theta^*,X_{\tau_i})(Y_{i+1}-\eta(\theta^*,X_{\tau_i})) \right)\qquad\text{for }t < 0}  \\ 
 \end{array} \right.
\end{equation}
\begin{theorem}\label{T:asiymptotic_normal_n}
  Let $Y$, and $Y_n$ be defined as in Theorem
  \ref{T:Normality_One_dimension}.  Assume the algorithm given by
  equation$~(\ref{E:NormalMultiRecursion})$ where $K$ is a nonsingular
  symmetric positive definite matrix. Let $\theta^*$ be an isolated
  stable point of the ODE $~(\ref{E:PODE})$ in the interior of $\mathbb{H}$,
  and assume that $(\varTheta_n)$ converges almost surely
  to the process with constant value $\theta^*$.  Assume that
  $\eta(\cdot,d)$, defined as above for $d\in\mathtt{D}$ are twice
  continuous differentiable functions.  Assume that the Hessian matrix
\[
A=D^2_{\theta^*}(\frac{1}{2}\sum_{d\in\mathtt{D}}
p_d(\eta(\theta,d)-\eta(\theta^*,d))^2)
\]
is positive definite.  Moreover, assume that the eigenvalues of the
matrix $KA$ are greater that $1/2$ (in particular the matrix
$(-KA+I/2)$ is negative definite).  Then, the sequence
$(U^n(\cdot),W^n(\cdot))$ converges weakly in
$D^r\left[0,\infty\right)\times D^r\left(-\infty,\infty\right)$ to a
limit process $\left(U(\cdot),W(\cdot)\right)$, where $W(\cdot)$ is a
Wiener process with covariance matrix
\[
\Sigma=\sum_{d\in\mathtt{D}}p_dVar_d(Y)(K\triangledown(\theta^*,d))(K\triangledown(\theta^*,d))^{\prime}
\]
and $U(\cdot)$ a stationary process with
\begin{eqnarray*}
U(t)&=U(0)+\int_{0}^{t}{(-KA+I/2)U(s)}\,ds +W(t)\\
&=\int_{-\infty}^{t}{\exp{((-KA+I/2)(t-s))}} \,dW(s)
\end{eqnarray*}
\end{theorem}
The proof of the previous theorem follows in a straightforward manner
from the rate of convergence theorem for stochastic algorithms with
exogenous noise and decreasing step size.  See for instance
\cite{hKgY97}, Theorem 10.2.2.  Details are left to the reader.  A few
words are needed to review the hypothesis of the previous theorem.
Theorem \ref{T:consistencyMulty} above, gives sufficient conditions
for the almost surely convergence of the sequence
$\left(\varTheta_n\right)$.  If $\theta^*$ is a unique stationary
point of the ODE $~(\ref{E:PODE})$ in the interior of $\mathbb{H}$, and $K$,
$A$ are positive definite symmetric matrices then $\theta^*$ is a
stable point of the ODE.  The latter follows from the basic theory of
Dynamical systems.  See for instance Perko \cite{Perko01}.  It is worth noting
that Theorem \ref{T:asiymptotic_normal_n} requires the existence of a
unique stationary point of the ODE in the interior of $\mathbb{H}$.  In the
setting of Theorem \ref{T:asiymptotic_normal_n}, a sufficient
condition for the uniqueness in the interior of $\mathbb{H}$ is discussed after Theorem \ref{T:consistencyMulty}.  For
any matrix $A$, as above, it is clearly possible to find a symmetric,
positive definite matrix ``large enough'' so that the eigenvalues of
the matrix are greater that $1/2$.

Among the choices for $K$ in equation$~(\ref{E:NormalMultiRecursion})$
the asymptotic optimal covariance is achieved by $K=A^{-1}$.  For
this, the limit $U(\cdot)$ satisfies
\[
dU=\left(-KA+I/2)\right)U\,dt + K\Sigma^{1/2}dW_0
\]
where $W_0$ is the standard Winner process.  The stationary covariance
is
\[
\int_{0}^{\infty}{e^{(-KA+I/2)t} K\Sigma
  K^{\prime}e^{(-A^{\prime}K^{\prime}+I/2)t}}\,dt
\]
the trace of this matrix is minimized by choosing $K=A^{-1}$, which
yields the asymptotically covariance $A^{-1}\Sigma (A^{\prime})^{-1}$.
See Kushner and Yin \cite{hKgY97} for a deeper discussion on the latter.  In order to
determine a choice of $K$ that is optimal for the class of estimators
proposed by equation$~( \ref{E:NormalMultiRecursion})$ it is necessary
to have a consistent estimator for the parameter $\theta^*$.  This can
be accomplished by initially employing a not necessarily optimal estimator
of the type mentioned above.  As we mentioned earlier, in section
\ref{S:Consistent}, it is possible to make use of Theorem
\ref{T:asiymptotic_normal_n} for the estimation of the parameters related
with the diffusion or the parameters related with the ratio between the
drift and the diffusion.  
\section{Examples}\label{S:examples}
In this section we show as an illustration the estimation for two
examples from the financial literature.  We choose the estimation for
the one-dimensional parameter space.  Details for the multidimensional
parameter space are left to the reader to fill out.

\subsection{CEV} \label{S:exonedim}
In this example we consider the estimation of the parameters for the
constant elasticity of variance (CEV) process introduced by
Cox \cite{jC75} and by Cox and Ross \cite{jCsR76}. The application of this process to
interest rates is discussed in Marsh and Rosenfeld \cite{tMeR83}.

Let us assume that
$(X_t,\mathcal{F}_t,\mathbf{P}_{x}^{(\lambda,\varsigma)})_{(\lambda,\varsigma)
  \in \mathbb{R}^2}$ is a parametric set of diffusions on $\mathbb{R}$
with sample space $(\Omega, \mathcal{F}_{\infty})$.  We assume that
the differential operator of the diffusion
$(X_t,\mathcal{F}_t,\mathbf{P}_x^{(\lambda,\varsigma)})$ is given by
the formula
\begin{equation*}
        L_{\lambda, \varsigma}= \frac{1}{2}
                \sigma^2(\varsigma) x^{2\gamma}\frac{d^2}{dx^2} +
        (\mu/\sigma^2)(\lambda) \sigma^2(\varsigma) \frac{d}{dx}
        \qquad \text{ for }  \lambda, \varsigma \in \mathbb{R}
\end{equation*}
where $\mu/\sigma^2\colon\mathbb{R}\mapsto\mathbb{R}^+$ is a
differentiable function such that $d (\mu/\sigma^2)/d \lambda >0$ and
$\sigma^2 \colon \mathbb{R} \mapsto \mathbb{R}^+$ is a continuous
differentiable function with $d (\sigma^2)/d \varsigma >0$ that
satisfies equation~(\ref{E:rate}), and $\gamma >1$ is a fixed known
constant.  Let $(\lambda^*,\varsigma^*)$ be the true parameter of the
process.  Let $\mathtt{D} =\{d_{1}, \dots ,d_{s}\} \subset \mathbb{R}$
be a finite set of real numbers where $d_1 < \dots <d_s$. First we
consider the collection of diffusions $(X_t,
\mathcal{F}_t,\mathbf{P}_{x}^{(\lambda,\varsigma^*)})_{\lambda \in
  \mathbb{R}}$.  Let $f \colon \mathbb{R} \to \mathbb{R}$ be the
identity function, and let $\eta^f, V^f$ be defined as in
equations~(\ref{E:average}) and~(\ref{E:deltaaverage}) respectively.
We suppress $f$ in what follows.  It follows that
\begin{equation*} 
        \eta(\lambda,x)
         =\{\mathtt{D}_r(x)-\mathtt{D}_l(x)\}
                \frac{
                  \int_{\mathtt{D}_l(x)}^x
                  \exp
                  ( \frac{\mu/\sigma^2(\lambda)}
                  {\gamma-1}
                   y^{2-2\gamma}) \ dy
                  }
                {
                  \int_{\mathtt{D}_l(x)}^{\mathtt{D}_r(x)}
                  \exp ( \frac{\mu/\sigma^2(\lambda)}{\gamma-1}
                  y^{2-2\gamma}
                  )\ dy
                  } + \mathtt{D}_l(x)
\end{equation*}
for $x \in \mathtt{D}$.  Let $(\Theta_{n})$ and $(\Theta^N_n)$ be the
$(\mathcal{F}_{\tau_n})$ adapted sequences of random variables defined
as in equations~(\ref{E:sequence}) and
~(\ref{E:Normal_Recursion_scale}) respectively, where $H \colon
\mathbb{R} \times \mathtt{D} \to \{1,-1\}$ is the constant function
taking value $-1$, $\alpha(d)=-(\partial \eta /\partial
\lambda)(\lambda^*,d)$ for $d \in \mathtt{D}$ and $\Theta_1$ is a
finite $\mathcal{F}_1$ measurable random variable.  We observe that
the computation of the estimators $\Theta_n$ does not depend on the
value of $\varsigma^*$.

It follows by Corollary~\ref{C:scalecons} that the sequence of
estimators $(\Theta_n)$ converges almost surely
$\mathbf{P}_x^{(\lambda^*, \varsigma^*)}$ to $\lambda^*$.  Let $\mu$
be a probability measure on $S$ and let $(\Omega, \mathbf{P},
\mathcal{F}_{\infty})$ be the sample space where $\mathbf{P}=\int
\mathbf{P}^{(\lambda^*,\varsigma^*)}_x d \mu(x)$.  As a consequence of
Corollary~\ref{C:scalenormal}, $n^{1/2}(\Theta^N_n -\lambda^*)$ is
asymptotically normally distributed with mean zero and variance
$\sigma^2 =\sum_{d \in \mathtt{D}}p_d
\mathbf{E}^{(\lambda^*,\varsigma^*)}_d (X_{\nu_2} -
\eta(\lambda^*,d))^2/\alpha^2(d)$ where $p=(p_i)$ is the left-fixed
probability row vector of the Markov chain $(X_{\tau_n},
\mathcal{F}_{\tau_n})$ as in Lemma~\ref{L:Ergodicity}.

Next, we assume that the parameter $\lambda^*$ is known. (The later
can be assumed by the previous remark.)  We consider the collection of
diffusions
$(X_t,\mathcal{F}_t,\mathbf{P}_{x}^{(\lambda^*,\varsigma)})_{\varsigma
  \in \mathbb{R}}$.  We define $\tilde{\eta}^g$ and $\tilde{V}^g$ by
equations~(\ref{E:timeaverage}) and ~(\ref{E:timedeltav}) respectively
where $g$ is the identity function.  See Lemma~\ref{L:avwaiting} of
Appendix B for a computation of these functions. We suppress $g$ in
what follows.  Let $(\tilde{\Theta}_n)$, $(\tilde{\Theta}^N_n)$ be the
$(\mathcal{F}_{\tau_n})$ adapted sequence of random variables as in
equations~(\ref{E:seqtime}) and ~(\ref{E:Normal_Recursion_timecons})
respectively, where $\tilde{\Theta}_1$ is a bounded $\mathcal{F}_1$
random variable, $\tilde{H} \colon \mathbb{R} \times \mathtt{D}
\mapsto \{-1,1\}$ is the constant function equal to 1 and
$\tilde{\alpha}(\cdot)=-(\partial\tilde{\eta}/\partial
\varsigma)(\varsigma^*,\cdot)$.  It follows by
Corollary~\ref{C:timecons} that for any $x \in \mathbb{R}$ the
sequence $(\tilde{\Theta}_n)$ converges almost surely
$\mathbf{P}^{(\lambda^*,\varsigma^*)}_x$ to $\varsigma^*$.  It follows
that the sequence of estimators $(\tilde{\Theta}_n)$ converges almost
surely $\mathbf{P}_x^{(\lambda^*, \varsigma^*)}$ to $\varsigma^*$.
Next, it follows by Corollary~\ref{C:timecons_normality}, that
$n^{1/2}(\tilde{\Theta}^N_n -\varsigma^*)$ is asymptotically normally
distributed with mean zero and variance $\sigma^2 =\sum_{d \in
  \mathtt{D}}p_d \mathbf{E}^{(\lambda^*,\varsigma^*)}_d (\nu_2 -
\tilde{\eta}(\varsigma^*,d))^2/\tilde{\alpha}^2(d)$.
\subsection{Cox-Ingersoll-Ross}  \label{S:examcox}
In this example we consider the estimation of some parameters for the
model of term structure of interest rates of Cox-Ingersoll-Ross.  See
Cox et al \cite{CIR85}.  We consider the problem of the estimation of the
quotient between the ``speed of adjustment'' and the ``volatility of
the process''. See Cox et al \cite{CIR85} for explanation of this terminology.
Let us assume that
$(X_t,\mathcal{F}_t,\mathbf{P}_{x}^{(\lambda,\varsigma)})_{(\lambda,\varsigma)
  \in \mathbb{R}^2}$ is a parametric set of diffusions on $\mathbb{R}$
with sample space $(\Omega, \mathcal{F}_{\infty})$. Let
$(\lambda^*,\varsigma^*)$ be the ``true'' parameter of the process. We
assume that the differential operator of the diffusion
$(X_t,\mathcal{F}_t,\mathbf{P}_x^{(\lambda,\varsigma)})$ is given by
the formula
\begin{equation*}
        L_{\lambda, \varsigma}= \frac{1}{2}
                \sigma^2(\varsigma)\, x\frac{d^2}{dx^2} +
        (\mu/\sigma^2)(\lambda) \sigma^2(\varsigma)(\alpha-x) \frac{d}{dx}
        \qquad \text{ for }  \lambda, \varsigma \in \mathbb{R}
\end{equation*}
where $\alpha \in \mathbb{R}^+$ is a given known constant,
$\mu/\sigma^2\colon\mathbb{R}\mapsto\mathbb{R}^+$ is a differentiable
function with $d (\mu/\sigma^2)/d \lambda >0$ and $\sigma^2 \colon
\mathbb{R} \mapsto \mathbb{R}^+$ is a continuous differentiable
function with $d (\sigma^2)/d \varsigma >0$ that satisfies equation
~(\ref{E:rate}). Let $\mathtt{D} =\{d_{1}, \dots ,d_{s}\} \subset
(0,\alpha) \cup (\alpha,\infty)$ be set of positive real numbers such
that $d_1< \dots <d_s$.  First, we consider the collection of
diffusions
$(X_t,\mathcal{F}_t,\mathbf{P}_x^{\lambda,\varsigma^*})_{\lambda\in\mathbb{R}}$.
Let $f \colon \mathbb{R} \to \mathbb{R}$ be the identity and let
$\eta^f, V^f$ be defined as in equations~(\ref{E:average})
and~(\ref{E:deltaaverage}) respectively.  We suppress $f$ in what
follows.  Moreover we assume that $U_d \subset [0, \alpha) \cup
(\alpha, \infty)$ for $d \in \mathtt{D}$.  Although a diffusion with
differential operator as above is not a regular Markov process on
$\mathbb{R}$, we observe that if $1/2\alpha > (\mu/\sigma^2)(\lambda)$
$(1/2\alpha \leq (\mu/\sigma^2)(\lambda))$ then the part of the
process on $[0,\infty)$ (on $(0,\infty)$) (see \cite{eD65}, volume I
for a definition of the part of a process on a subset of the state
space) is a regular diffusion on $[0,\infty)$ (on $(0,\infty)$) in the
sense of definition 15.1 of  Dynkin \cite[volume II, p.121]{eD65}.  For a
discussion of this see for example Cox et al \cite{CIR85}.  Either case follows
from the analysis of the boundary classification criteria, see for
instance Gihman and A.V. Skorohod \cite{GS72}.  In either case Condition~\ref{C:regularity} is
satisfied if $a,b \in (0,\infty)$.  It follows that
\begin{equation*} 
        \eta^f (\lambda,x)
         =\{\mathtt{D}_r(x)-\mathtt{D}_l(x)\}
                \frac{\int_{\mathtt{D}_l(x)}^x y^{-2\alpha(\mu/
                \sigma^2)(\lambda)} \exp(2(\mu/
                \sigma^2)(\lambda)y) \, dy}
                {\int_{\mathtt{D}_l(x)}^{\mathtt{D}_r(x)} y^{-2\alpha(\mu/
                \sigma^2)(\lambda)} \exp(2 (\mu/
                \sigma^2)(\lambda)y)\, dy})  + \mathtt{D}_l(x)
\end{equation*}
for $x \in \mathtt{D}$.  Let $\Theta_1$ be a finite $\mathcal{F}_1$
random variable, and let $(\Theta_{n})$ be the
$(\mathcal{F}_{\tau_n})$ adapted sequence of random variables defined
as in equation~(\ref{E:sequence}) where $H \colon \mathbb{R} \times
\mathtt{D} \to \{1,-1\}$ is the function defined by the formula
$H(\lambda,d)=\mathbf{1}_{(-\infty,\alpha)}(d)-
\mathbf{1}_{(\alpha,\infty)}(d)$. If $\lambda^* \in \mathbb{R}$ is a
fixed number it follows by Corollary~\ref{C:scalecons} that the
sequence $(\Theta_n)$ of random variables converges to $\lambda^*$
almost surely $\mathbf{P}_{x}^{(\lambda^*,\varsigma^*)}$ for any $x
\in (0,\infty)$.
Next, we assume that the parameter $\lambda^*$ is known.  We consider
the collection of diffusions $(X_t,\mathcal{F}_t,\mathbf{P}_x^{(\lambda^*,\varsigma)})_{\varsigma\in\mathbb{R}}$. We define $\tilde{\eta}^g$ and $\tilde{V}^g$ by
equations~(\ref{E:timeaverage}) and ~(\ref{E:timedeltav})
respectively, where $g$ is the identity function, $\mathtt{D} \subset
(0,\infty)$ and $U_d \subset (0,\infty)$.  Let $(\tilde{\Theta}_n)$ be
a sequence of random variables as in equation~(\ref{E:seqtime}) where
$\tilde{\Theta}_1$ is a bounded random variable and $\tilde{H} \colon
\mathbb{R} \times \mathtt{D} \mapsto \{-1,1\}$ is the constant
function equal to 1.  It follows by Corollary~\ref{C:timecons} that
for any $\varsigma^* \in \mathbb{R}$, the sequence
$(\tilde{\Theta}_n)$ converges almost surely
$\mathbf{P}^{(\lambda^*,\varsigma^*)}_x$ to $\varsigma^*$ for any $x \in(0,\infty)$.

Similar considerations can be made if we want to estimate the
``central location'' or ``long term value'' of the process and the
diffusion coefficient.  See Cox et al \cite{CIR85} for a explanation of this terminology.  

Last, we notice that as a consequence of Corollary~\ref{C:scalenormal}
and Corollary~\ref{C:timecons_normality} asymptotic normality of the
appropriate normalization of the estimators constructed for the
Cox-Ingersoll-Ross model can be obtained.  The details are left to the
reader.

\section{Conclusion} \label{S:Conclusion}
The thrust of this paper has been to introduce the ideas of stochastic
algorithms to the problem of the estimation of parameters of a
continuous diffusion process using observed discrete data.  The later
could be potentially useful for the study of non time-homogeneous
diffusion processes.   Besides, we have proposed  sampling schemes
that depends on space discretization rather than time discretization.
These sampling schemes are closer to the Markov character of diffusion
processes.  We also propose a new parameterization of diffusions that we
believe is closer in spirit to the
initial attempts made in probability to describe a diffusion by its
``road map'' and ``speed''.  The main
results given here (construction of sequences of estimators,
asymptotic consistency of such sequences, and asymptotic normality of
such sequences) as well as the two examples taken from Mathematical
Finance dealt with families of diffusion processes that have a
one-dimensional state space and a multidimensional parameter space.

Future questions will center on the generalization of the current
techniques for use in the case of a multi-dimensional state space and
the development for the current setting of stochastic algorithms
appropriate to the description of non time-homogeneous
diffusions. A 
particularly interesting question is to find sufficient conditions on
multi-dimensional parameterizations of diffusion operators that
guarantee identification of the corresponding process from moment
conditions of the type presented in this paper. Another direction of
research can be centered on the the effective computation of
the stochastic algorithms presented here and its comparison with
other techniques.
\section*{Acknowledgement}

I am very grateful to  my advisor Professor Neil E. Gretsky for his
patience and constructive suggestions.  I also want to thank
Professor Darrel Duffie for his suggestion about my area of
research.  All remaining errors are mine.

\providecommand{\bysame}{\leavevmode\hbox to3em{\hrulefill}\thinspace}
\providecommand{\MR}{\relax\ifhmode\unskip\space\fi MR }
\providecommand{\MRhref}[2]{%
  \href{http://www.ams.org/mathscinet-getitem?mr=#1}{#2}
}
\providecommand{\href}[2]{#2}

\appendix
\section{Appendix. Proofs}
\begin{proof}[Proof of Theorem \ref{T:Monro}]
  If we define $T_n=\varTheta_n-\theta^*$ then equation
  ~(\ref{E:Recursion}) becomes
\begin{equation*}
  T_{n+1}=T_n-\gamma_n H(\varTheta_n, Y_n)V(\varTheta_n, Y_n,X_{n+1})
\end{equation*}
It follows that
\begin{multline} \label{E:proof2}
  \mathbf{E}(\|T_{n+1} \|^2 \mid \mathcal{F}_n) - \|T_n\|^2=\ -2
  \gamma_n T_n \cdot H(\varTheta_n,Y_n) \overline{V}
  (\varTheta_n,Y_n)\\
  + \gamma^2_n S^2 (\varTheta_n, Y_n) \leq \gamma^2_n K (1 + \|T_n \|
  ^2)
\end{multline}
where the last inequality follows by equations~(\ref{E:H1a}) and
~(\ref {E:H2}).  Moving the terms that have either $\|T_n\|$ or
$\|T_{n+1}\|$ to the left of the previous equation we obtain
\begin{equation} \label{E:proof3}
  \mathbf{E}(\|T_{n+1} \|^2 \mid \mathcal{F}_n) -\|T_n\|^2 (1 + K
  \gamma^2_n) \leq K \gamma^2_n
\end{equation}
Define $\Pi_n=\prod_{k=1}^{n-1}(1+K \gamma^2_n)$ and let
$T_{n}'=(\frac{1}{\Pi_n})^ {\frac{1}{2}} T_n$.  We observe that the
sequence $(\Pi_n)$ is a convergent sequence of positive numbers since
$(\log \Pi_n)$ converges by Hypothesis $H_3$.  Using equation~(\ref{E:proof3}) 
we obtain
\begin{equation} \label{E:proof4}
  \mathbf{E}(\|T_{n+1} ' \|^2 \mid \mathcal{F}_n) -\|T_n' \|^2 \leq K
  \frac{\gamma^2_n}{\Pi_{n+1}}
\end{equation}
If $F_n=\{ \omega \in \Omega \mid \mathbf{E}( \| T_{n+1}' \|^2 \mid
\mathcal{F}_n) - \|T_n' \|^2 \geq 0 \}$, then equation~(\ref{E:proof3}) 
and Hypothesis $H_3$ imply that
\begin{equation*}
  \sum_{n=1}^{\infty} \mathbf{E} (\mathbf{1}_{F_n} (\|T_{n+1} ' \|^2-
  \|T_n ' \|^2)) < \infty
\end{equation*}
It follows by the almost sure convergence of quasi-martingales  that $T_n^{'}$
converges almost surely toward a positive integrable random variable (see, for example, Theorem 9.4
  page 49 and Proposition 9.5 of  M\'{e}tivier \cite{mM82}).
We conclude that the same property holds for $T_n$.  The next step of
the proof is to prove that the convergence of $T_n$ is to zero.  By
inequality~(\ref{E:proof4}), Hypothesis $H_3$, the definition of
$T_n$ and the fact that $\varTheta_1$ belongs to $L^2(\mathbf{P})$, it
follows that
\begin{equation} \label{E:proof6}
  \sup_n \mathbf{E}( \|T_n\|^2) < \infty
\end{equation}
We also observe that
\begin{align*}
  0 & \leq\sum_{n=1}^{\infty} 2 \gamma_n \mathbf{E} (T_n \cdot
  H(\varTheta_n,Y_n)
  \overline{V}(\varTheta_n,Y_n))  \\
  & \leq \sum_{n=1}^{\infty} \mathbf{E}( \| T_n \|^2) -
  \mathbf{E}(\|T_{n+1}\|^2) + (1 + \sup_{k \geq 0}
  \mathbf{E}(\|T_k\|^2)) \sum_{n=1}^{\infty} K \gamma_n^2 < \infty
\end{align*}
where the last inequality follows by equations~(\ref{E:proof2}), ~(\ref{E:proof6}), and Hypothesis $H_3$.  Since
$\sum_{k} \gamma_{n_k}= \infty$ there exists a subsequence $(n_k')$ of
$(n_k)$ such that
\begin{equation} \label{E:proof9}
  \lim_{k} \mathbf{E}(T_{n_k'} \cdot H(\varTheta_{n_k'},Y_{n_k'})
  \overline{V}(\varTheta_{n_k'}, Y_{n_k'})=0
\end{equation}
equation~(\ref{E:proof9}) implies that for any $\varepsilon >0$,
$\liminf_{k} \|T_{n_k^{'}} \| \leq \varepsilon$ almost surely.  To
prove this last statement, let us assume otherwise.  Then there exists
$\varepsilon >0$ such that $T_{n'_k} \geq \varepsilon$ for all $k $
big enough on some set $A$ of probability greater than zero.  It would
follow by Fubini's theorem and equation~(\ref{E:H1b}) that there
exists $\delta >0$ such that
\begin{multline*}
  \int\limits_{\Omega} T_{n'_k} \cdot
  H(\varTheta_{n'_k},Y_{n'_k})\overline{V}
  (\varTheta_{n'_k},Y_{n'_k})  d \mathbf{P}\\
  \geq \int\limits_{A} T_{n'_k} \cdot H(\varTheta_{n'_k}, Y_{n'_k})
  \overline{V}(\varTheta_{n'_k},Y_{n'_k})\, d \mathbf{P} \geq
  \int\limits_{A} \delta d \mathbf{P} \geq \delta \mathbf{P}(A) >0
\end{multline*}
for all $k $ big enough.  The former is in contradiction with equation~(\ref{E:proof9}).  
Since $(T_n)$ converges almost surely, then $T_n
\to 0$ almost surely.
\end{proof}
\begin{proof}[Proof of Corollary \ref{C:scalecons}]
  Let $\mu$ be a probability measure on $\mathbb{R}$ supported on $S$,
  and let $\mathbf{P}=\int \mathbf{P}^{\lambda^*}_x \, d \mu$.  We
  define $\overline{V} \colon \mathbb{R} \times \mathtt{D} \mapsto
  \mathbb{R}$ and $S^2 \colon \mathbb{R} \times \mathtt{D} \mapsto
  \mathbb{R}$ by the following formulas:
\begin{gather}
  \overline{V}(\lambda, d)= \eta^{f}(\lambda^*,d)-
  \eta^{f}(\lambda,d) \label{E:cons1}\\
  S^2(\lambda,d)=(\eta^{f}(\lambda^*,d)- \eta^f(\lambda,d))^2
  +(\eta^{f^2}(\lambda^*,d) - (\eta^f(\lambda^*,d))^2)
\end{gather}
It follows by the strong Markov property of
$(X_t$,$\mathcal{F}_t$,$\mathbf{P}_{x}^{\lambda^*})$ that Conditions
$A_1$ and $A_2$ of Theorem~\ref{T:Monro} are satisfied. Since
$\eta^f(\cdot,d)$ for any $d \in \mathtt{D}$ is bounded and
$\mathtt{D}$ is finite it follows that Property $H_2$ of
theorem~\ref{T:Monro} is satisfied.  By Corollary~\ref{C:criteria} of
Appendix B, equation~(\ref{E:H1a}) of Theorem~\ref{T:Monro} is
satisfied.  Last, we notice that \small
\begin{equation} 
  \mathbf{E} ((\lambda-\lambda^*)
  H(\lambda,X_{\tau_n})\overline{V}(\lambda,X_{\tau_n})) =
  \sum_{m=1}^{s} (\lambda-\lambda^*) \overline{V}(\lambda,d_m)
  \mathbf{P}(X_{\tau_n}=d_m)
\end{equation}
\normalsize By the last equation and Corollary~\ref{C:appendixA1} of
Appendix A, equation~(\ref{E:H1b}) holds.  It follows by
Theorem~\ref{T:Monro} that the sequence of random variables
$(\varTheta_n)$ converges almost surely $\mathbf{P}$ to $\lambda^*$.
Since the last statement holds for any initial probability measure
supported on $S$ the result follows.
\end{proof}
\begin{proof}[Proof of Corollary \ref{C:timecons}]
  Let $\mu$ be a probability measure on $\mathbb{R}$ supported on $S$,
  and let $\mathbf{P}=\int \mathbf{P}^{\varsigma^*}_x \, d \mu$.  If
  we define $\widetilde{\overline{V}} \colon \mathbb{R} \times
  \mathtt{D} \mapsto \mathbb{R}$ and $\widetilde{S}^2 \colon
  \mathbb{R} \times \mathtt{D} \mapsto \mathbb{R}$ by:
\begin{gather}
  \widetilde{\overline{V}}(\varsigma, d)=
  \widetilde{\eta}^{f}(\varsigma^*,d)-
  \widetilde{\eta}^{f}(\varsigma,d) \label{E:timecons1}\\
  \widetilde{S}^2(\varsigma,d)=(\widetilde{\eta}^{f}(\varsigma^*,d)-
  \widetilde{\eta}^f(\varsigma,d))^2
  +(\widetilde{\eta}^{f^2}(\varsigma^*,x) -
  (\widetilde{\eta}^f(\varsigma^*,x))^2)
\end{gather}
It follows by the strong Markov property of
$(X_t,\mathcal{F}_t,\mathbf{P}_{x}^{\varsigma^*})$ that Conditions
$A_1$ and $A_2$ of Theorem~\ref{T:Monro} are satisfied.  By
Assumptions~\ref{P:scale} and \ref{P:separability} and
Lemma~\ref{L:avwaiting} of Appendix B, it follows that Property $H_2$
of Theorem~\ref{T:Monro} is satisfied.  The monotonicity and
differentiability of $\sigma_0$, Property~\ref{P:separability} of
Corollary~\ref{C:timecons} and Lemma~\ref{L:avwaiting} of Appendix B
imply equation~(\ref{E:H1a}) of Theorem~\ref{T:Monro}.  Last, we
notice that
\begin{equation} 
  \mathbf{E} ((\varsigma-\varsigma^*)
  \widetilde{H}(\varsigma,X_{\tau_n})\widetilde{\overline{V}}
  (\varsigma,X_{\tau_n})) = \sum_{m=1}^{s} (\varsigma-\varsigma^*)
  \widetilde{\overline{V}}(\varsigma,d_m) \mathbf{P}(X_{\tau_n}=d_m)
\end{equation}
By the last equation and Corollary~\ref{C:appendixA1} of Appendix A,
equation~(\ref{E:H1b}) holds.  It follows by Theorem~\ref{T:Monro}
that the sequence of random variables $(\widetilde{\Theta}_n)$
converges almost surely $\mathbf{P}$ to $\varsigma^*$.  Since the last
statement holds for any initial probability measure supported on $S$
the result follows.
\end{proof}

  The following lemma is used in the proof of Theorem
  \ref{T:consistencyMulty}, Theorem \ref{T:asiymptotic_normal_n} and
  Theorem  \ref{T:Martingale}. 
\begin{lemma} \label{L:Ergodicity} 
  Let $f \colon \mathtt{D} \to \mathbb{R}$ and $A=(a_{i,j})$ be the
  irreducible transition matrix of the Markov chain $(X_{\tau_n},
  \mathcal{F}_{\tau_{n}})$ with left-fixed probability row vector
  $p=(p_i)>>0$.
        
  Then for any $x \in S$
\begin{equation} \label{E:Normality2.1}
        n^{-1} \sum_{k=1}^{n} f(X_{\tau_k})  \to \sum_{d \in \mathtt{D}}
                f(d)p_d  \qquad \text{a.s. } \mathbf{P}_x
\end{equation}
\end{lemma}
\begin{remark}
  $A$ is a irreducible matrix. This follows by the recurrence of the
  Markov chain $(X_{\tau_n}, \mathcal{F}_{\tau_n})$.  Hence, there
  exists a unique left-fixed probability vector $p$; moreover, the
  entries of $p$ are strictly positive.  See for instance
   Petersen \cite{kP83}.
\end{remark}
\begin{proof}
  Let $\tilde{\mathbf{P}}_d$ for $d \in \mathtt{D}$ be the restriction
  of $\mathbf{P}_d$ to $\sigma( X_{\tau_1}, X_{\tau_2}, \cdots)$ and
  let $\mathbf{P}$ be the probability measure defined on $\sigma(
  X_{\tau_1}, X_{\tau_2}, \cdots)$ by the formula $\mathbf{P} =
  \sum_{d \in \mathtt{D}} p_d \tilde{\mathbf{P}}_d$.  By the strong
  Markov property $\theta_{\tau_2}$ defines a measure-preserving
  transformation.  Indeed, $(X_{\tau_n})$ is an irreducible Markov
  shift.  It follows by the point-wise ergodic theorem and the fact
  that an irreducible Markov shift is ergodic that equation~(\ref{E:Normality2.1}) 
holds a.s. $\mathbf{P}$ (See
   Petersen \cite{kP83}); the result follows using the strong Markov
  property and the fact that the components of the invariant vector
  $p$ are positive.
\end{proof}

\begin{theorem} \label{T:Martingale}
  Let $Z\colon (\Omega, \mathcal{F}_{\tau_2}) \to (\mathbb{R},
  \mathcal{B}(\mathbb{R}))$ be a measurable map that is bounded below.
  Moreover, assume that $Z \in \bigcap_{d \in \mathtt{D}}
  L^2(\mathbf{P}_d)$ and $\mathbf{E}_d (Z)=0$ for $d \in \mathtt{D}$.
  Let $Z_n$ be defined as $Z_n=Z \circ \theta_{\tau_n}$, for $n \geq
  2$ and $Z_1=Z$.  Then for any initial probability measure, the
  distribution of
\begin{equation} \label{E:Martingale1}
        n^{-1/2} \sum_{k=1}^n Z_n
\end{equation}
approaches the normal distribution with mean zero  and variance $\sigma^2=$\\
$\sum_{d \in \mathtt{D}} p_d \mathbf{E}_d(Z^2)$ where $p=(p_i)$ is the
left-fixed probability row vector of the Markov chain $(X_{\tau_n},
\mathcal{F}_{\tau_{n}})$ as in Lemma~\ref{L:Ergodicity}.
\end{theorem}
\begin{proof}
  Let $\mu$ be a probability measure on $S$, and let
  $(X_t,\mathcal{F}_t,\mathbf{P}_{\mu})$ be the Markov process with
  initial probability measure $\mu$.  We observe that $(Z_n)$ is a
  $(\mathcal{F}_{\tau_{n+1}})$ adapted process and
        \begin{gather}
          \mathbf{E}_{\mu}(Z_n \mid \mathcal{F}_{\tau_n})=0
                \label{E:Martingale2}\\
                \mathbf{E}_{\mu}(Z^2_n)=
                \mathbf{E}_{\mu}\mathbf{E}_{X_{\tau_n}}(Z^2) < \infty
                \label{E:Martingalel2}
        \end{gather}
        for $n \geq 1$.  The proof of the theorem follows using the
        strong Markov property, Lemma~\ref{L:Ergodicity}, and a line
        of argument similar to the technique of the proof of the
        Lindeberg-L\'{e}vy theorem for martingales (See
         Billingsley \cite{pB61}).
\end{proof} 
For the proof of Theorem~\ref{T:Normality_One_dimension} we make use
of the following easily proved lemma.

\begin{lemma} \label{L:Normality1}
  Let $\{\beta_{k,n} \mid 0 \leq k \leq n\}$ be the double indexed
  sequence of positive numbers defined as $\beta_{k,n}=\prod_{i=k+1}^n
  (1-1/i)$.  Then
\begin{equation} \label{E:Normality2}
        (1-\epsilon_k) \frac{k}{n} \leq \beta_{k,n} \leq (1+ \epsilon_k )\frac{k}{n}
\end{equation}
where $(\epsilon_k)$ is a sequence of positive numbers such that
$\epsilon_k \to 0$ as $k \to \infty$.  In particular $n^{1/2}
\beta_{k,n} \to 0$ as $n \to \infty$ for any fixed $k$.
\end{lemma}

\begin{proof}[Proof of Theorem \ref{T:Normality_One_dimension}]
  We observe that the sequence of random variables $(\Theta_n^N)$
  converges almost surely $\mathbf{P}$ to $\theta^*$ by Theorem~\ref{T:Monro}.
  Let $\overline{V} \colon \mathbb{R} \times
  \mathtt{D} \to \mathbb{R}$ and $S^2 \colon \mathbb{R} \times
  \mathtt{D} \mapsto \mathbb{R}$ be defined by the formulas:
\begin{gather*}
  \overline V(\theta, d)= \eta(\theta^*,d)-
  \eta(\theta,d)\\
  S^2(\theta,d)=(\eta(\theta^*,d)- \eta(\theta,d))^2 +
  \mathbf{E}_d(Y-\eta(\theta^*,d))^2
\end{gather*}
By Hypothesis $N_1$ it follows that for $\theta \in \mathbb{R}$, $d
\in \mathtt{D}$,
\begin{equation*}
(\theta- \theta^*) \frac{\overline{V}(\theta,d)}
{(\partial\overline{V}/
\partial \theta)(\theta^*,d)} \geq 0 
\end{equation*}
Hypothesis $N_2$ and the fact that $\alpha$ is defined on a finite set
and is nowhere zero imply that $S^2/\alpha^2$ satisfies Hypotheses
$H_2$ of Theorem~\ref{T:Monro}.  Using the strong Markov property
and an argument similar to the one used in the proof of Corollary~\ref{C:scalecons} we can prove that $(\Theta_n^N)$ converges almost
surely to $\theta^*$.  Let $Z \colon \mathtt{D} \times \mathbb{R} \to
\mathbb{R}$ be defined as $Z(d,y)=y-\eta(\theta^*,d)$.  We observe
that $Z(d,y)=V(\theta,d,y)-\overline{V}(\theta,d)$ for any $\theta \in
\mathbb{R}$.  We denote as $\delta \colon \mathbb{R} \times \mathtt{D}
\mapsto \mathbb{R}$ the function defined by the formula
$\overline{V}(\theta,d)=(\partial /\partial \theta)
\overline{V}(\theta^*,d)(\theta-\theta^*) +\delta(\theta-
\theta^*,d)$.  If we define $T^N_n=\Theta^N_n - \theta^*$ for $n \geq
1$, it follows that
\begin{equation} \label{E:Normality3}
        T^N_{n+1}=(1-\frac{1}{n}) T^N_n -
                \frac{Z_n}{n \alpha_n}
                -\frac{\delta_n}{n \alpha_n}
\end{equation}
where $Z_n=Z(X_{\tau_n},Y_{n+1})$, $\delta_n=\delta(T^N_n,X_{\tau_n})$
and $\alpha_n=\alpha(X_{\tau_n})$.  Iteration of equation~(\ref{E:Normality3}) yields
\begin{equation} \label{E:Normality4}
        T^N_{n+1}=\beta_{0,n} T^N_1 - \sum_{k=1}^n \frac{\beta_{k,n}}{k}
                 \frac{\delta_k}{\alpha_k}-
                \sum_{k=1}^n \frac{\beta_{k,n}}{k} \frac{Z_k}{\alpha_k}
\end{equation}
Hence, we can prove that $n^{1/2} T^N_n$ is asymptotically normal with
mean zero and variance $\sigma^2$, by proving that
\begin{gather}
  n^{1/2} \beta_{0,n} T_1^N \to 0 \quad \text{ almost surely,}
  \tag{ \ref{E:Normality4}a}\\
  n^{1/2} \sum_{k=1}^n \frac{\beta_{k,n}}{k} \frac{\delta_k}{\alpha_k}
  \to
  0 \quad \text{ in probability,} \tag{\ref{E:Normality4}b}\\
  n^{1/2} \sum_{k=1}^n (\frac{\beta_{k,n}}{k}-\frac{1}{n})
  (\frac{Z(X_{\tau_k},Y_{k+1})}{\alpha_k}) \to 0
  \quad \text{ in probability,} \tag{\ref{E:Normality4}c}\\
  n^{-1/2} \sum_{k=1}^n \frac{Z(X_{\tau_k},Y_{k+1})}{\alpha_k} \to
  N(0,\sigma^2) \quad \text { in distribution.}
  \tag{\ref{E:Normality4}d}
\end{gather}
Equation~(\ref{E:Normality4}a) follows by Lemma~\ref{L:Normality1}.  
Next, we observe that the terms on the left
hand side of equation~(\ref{E:Normality4}c) are uncorrelated.  We
observe that $ Y \in \bigcap_{d \in \mathtt{D}} L^2(\mathbf{P}_d)$ and
$\alpha$ is a nowhere zero function defined on a finite set. It
follows by the strong Markov property, and Lemma~\ref{L:Normality1}
that there exists a constant $C>0$ such that
\begin{equation} \label{E:Normality7}
        \mathbf{E}(n^{1/2} \sum_{k=1}^n (\frac{\beta_{k,n}}
        {k}-\frac{1}{n})
        (\frac{Z(X_{\tau_k},Y_{k+1})}{\alpha(X_{\tau_k})}))^2 \leq
        \frac{C}{n} \sum_{k=1}^{n} \epsilon^2_k
\end{equation}
where the right-hand side of equation~(\ref{E:Normality7}) goes to
zero as $n \to \infty$ since $\epsilon_k \to 0$ as $k \to \infty$.
Convergence in equation~(\ref{E:Normality4}c) follows by
Chebyshev's inequality.
Next we prove the convergence of equation~(\ref{E:Normality4}b).
We observe that equation~(\ref{E:proof4}) and the proof of Theorem~\ref{T:Monro} imply
\begin{equation} \label{E:Normality8}
        \limsup_{n} n \,\mathbf{E} (T^N_{n})^2 < \infty
\end{equation}
Let $\epsilon_1, \epsilon_2 >0$, and let $\overline{\delta}(\cdot)=
\max_{d \in \mathtt{D}} (\delta(\cdot,d)/\alpha(d))$.  Since
$\overline{\delta}(x)= o(\mid x \mid)$ there exists $\epsilon' >0$
such that
\begin{equation} \label{E:Normality9}
        \mid \overline{\delta}(x)
        \mid \leq \epsilon^2_2 \mid x \mid \qquad
        \text{ for } \mid x \mid \leq \epsilon' 
\end{equation}
Since $T^N_n \to 0$ almost surely, there exist $N_1 >0$ such that
\begin{equation} \label{E:Normality10}
        \mathbf{P}( \mid T^N_k \mid \leq \epsilon', k \geq N_1) >
        1-\epsilon_1
\end{equation}
It follows using equation~(\ref{E:Normality10}), the triangle
inequality, equation~(\ref{E:Normality9}), Markov's inequality, and
Lyapounov's inequality that
\begin{align} \label{E:Normality11}
  & \mathbf{P}(\mid n^{1/2} \sum_{k=N_1}^n \frac{\beta_{k,n}}{k}
  \frac{\delta_k}{\alpha_k} \mid > \epsilon_2) & \notag\\
  & \leq \epsilon_1 + \mathbf{P} (\mid n^{1/2} \sum_{k=N_1}^n
  \frac{\beta_{k,n}}{k} \frac{\delta_k}{\alpha_k} \mid > \epsilon_2,
  \mid T^N_k \mid \leq \epsilon', k \geq N_1) \notag\\
  & \leq \epsilon_1 + \mathbf{P}( n^{1/2} \sum_{k=N_1}^n
  \frac{\beta_{k,n}}{k} \mid \frac{\delta_k}{\alpha_k} \mid >
  \epsilon_2
  \mid,  T^N_k \mid \leq \epsilon', k \geq N_1) \notag \\
  & \leq \epsilon_1 + \mathbf{P}( \epsilon^2_2 n^{1/2} \sum_{k=N_1}^n
  \frac{\beta_{k,n}}{k} \mid T^N_k \mid
  > \epsilon_2)\\
  & \leq \epsilon_1 + \epsilon_2 \,\mathbf{E} (n^{1/2} \sum_{k=N_1}^n
  \frac{\beta_{k,n}}{k} \mid T^N_k \mid) \notag\\
  & \leq \epsilon_1 + \epsilon_2 ( \sum_{k=N_1}^n
  \frac{\beta_{k,n}}{k} n^{1/2} \mathbf{E}^{1/2} (T^N_k)^2 ) \notag
\end{align}
The convergence in probability in equation~(\ref{E:Normality4}b)
follows using equation~(\ref{E:Normality11}) and Lemma~\ref{L:Normality1}.  
Finally, we observe that the convergence in
distribution of equation~(\ref{E:Normality4}d) is a consequence of
Theorem~\ref{T:Martingale}.
\end{proof}
\begin{proof}[Proof of Corollary \ref{C:scalenormal}]
  We observe that $X_{\nu_n}=X_{\nu_2} \circ \theta_{\tau_{n-1}}
  =X_{\nu_1} \circ \theta_{\tau_2}^{n-2}$ for $n \geq 3$.  Indeed,
  $\eta$ satisfies Condition $N_1$ of Theorem~\ref{T:Normality_One_dimension} by Corollary~\ref{C:criteria},
  and by the definition of $\eta$ Condition $N_2$ is satisfied.  The
  result is a straightforward consequence of Theorem~\ref{T:Normality_One_dimension}.
\end{proof}
\begin{proof}[Proof of Corollary \ref{C:timecons_normality}]
  We observe that $\nu_n-\tau_{n-1}=\nu_2 \circ \theta_{\tau_{n-1}} =\nu_2 \circ
  \theta^{n-2}_{\tau_{2}}$ for $n \geq 3$.  By Lemma~\ref{L:avwaiting}
  $\tilde{\eta}$ satisfies Condition $N_1$ of
  Theorem~\ref{T:Normality_One_dimension}, and $\tilde{\eta}$
  satisfies Condition $N_2$ by Lemma~\ref{L:avwaiting}.  The result is
  a straightforward consequence of
  Theorem~\ref{T:Normality_One_dimension}.
\end{proof}

\section{Appendix} \label{S:Matrices}
In this appendix we derive some technical results about the
transition matrices of the Markov chain $(X_{\tau_n},
\mathcal{F}_{\tau_n})$. In this section all the matrices are
stochastic matrices. We state some easily proved
results.  The proof is left to the reader.\\

\begin{definition}
We say that a matrix $A$ of size $s \times s$ is of type
$\mathbf{I}$  if for all $i, \, j \in \{1 \dots s\}$, $i \equiv j
\mod{s}$ implies $a_{i,j}=0$ .   We would say that a matrix $A$ of
size $s \times s$ is of type $\mathbf{II}$ if whenever $i \equiv j
+1 \mod{s}$ implies $a_{i,j}=0$
for all  $i, \, j \in \{1 \dots s\}$\\
\end{definition}

\begin{lemma} \label{l:appendix1}
let $A$ and $B$ be two $s \times s$ matrices. Then the following holds: \\
If $A$ and $B$ are both matrices of type $\mathbf{I}$ then $AB$
is a matrix of type $\mathbf{II}$. \\
If $A$ and $B$ are matrices of type $\mathbf{II}$ then so is
$AB$.\\
If $A$ is a matrix of type $\mathbf{I}$ and $B$ is a matrix of
type $\mathbf{II}$ then $AB$ and $BA$ are  matrices of type
$\mathbf{I}$.\\
\end{lemma}

\begin{definition}
Given a matrix $A$ of type $\mathbf{II}$ we define $P_{even}(A)$
and $P_{odd}(A)$ to be the matrices of size $[\frac{s}{2}] \times
[ \frac{s}{2}]$, and size $[\frac{s+1}{2} ] \times [ \frac{s+1}{2}
]$ respectively defined by the following formulas
\begin{alignat*}{2} 
  (P_{even}(A))_{i,j}=a_{2i,2j} &&\quad \text{for } i,j
  \in \{1 \dots [\frac{s}{2}] \}\\
  (P_{odd}(A))_{i,j}=a_{2i-1,2j-1} &&\quad \text{for } i,j \in \{1
  \dots [\frac{s+1}{2}] \}\notag
\end{alignat*}
\end{definition}

\begin{lemma} \label{l:appendix2}
If $A$ and $B$ are $s \times s$ matrices of type $\mathbf{II}$
then the following property holds:
\begin{align*}
P_{even}(A)P_{even}(B)=P_{even}(AB)\\ 
P_{odd}(A)P_{odd}(A)=P_{odd}(AB) 
\end{align*}
In particular for any $n$ positive integer
\begin{align*}
P_{even}(A^n)=(P_{even}(A))^n \\
P_{odd}(A^n)=(P_{odd}(A))^n 
\end{align*}
\end{lemma}

It is obvious that a matrix $A$ of type $\mathbf{II}$ is
completely
determined by $P_{odd}(A)$ and $P_{even}(A)$.\\

\begin{lemma} \label{l:appendix3}
Let $A=(a_{i,j})$ be an $s \times s$ matrix whose entries are
non-negative and such that $a_{i,j} \ne 0$ for $\mid i-j \mid \leq
1$. Then for any positive integer $n$, $A^n$ is a matrix such that
$a^n_{i,j} \ne 0$ for $\mid i -j\mid \leq n$ where
$A^n=(a^n_{i,j})_{i,j}$
\end{lemma}

We observe that if $A$ is the transition  matrix of the Markov
process
 $(X_{\tau_n}, \mathcal{F}_{\tau_n})$ then $P_{odd}(A^2)$ and
$P_{even}(A^2)$ satisfies the condition of
the previous lemma.

\begin{corollary} \label{C:appendixA1}
If A is an $s \times s$ transition matrix of a Markov process
$(X_{\tau_n},\mathcal{F}_{\tau_n})$, then $P_{odd}(A^{2n})$ and
$P_{even}(A^{2n})$ converge to stochastic matrices $A_1$ and
$A_2$.  Indeed, there are matrices $C_1$ and $C_2$ where $C_1$ is
a $1 \times [ \frac{s+1}{2} ]$ matrix and $C_2$ is a $1 \times [
\frac{s}{2} ]$ matrix  such that  $A_1=(1,\dots,1)'C_1$ and
$A_2=(1, \dots, 1)'C_2$ and the components of $C_1$ and $C_2$ are
positive.
\end{corollary}
\begin{proof}
The result follows from Lemma~\ref{l:appendix2},
Lemma~\ref{l:appendix3} , the observation made right after the
proof of Lemma~\ref{l:appendix3} and the fundamental theorem for
regular Markov chains.  See for example  Kemeny and Snell \cite{jKjS76} Theorem
4.1.4.
\end{proof}
\section{Appendix} \label{S:Solution1}
In this appendix we state and prove some results that are needed for
the actual computation of the moments required for the construction of
the algorithms proposed.  Let us assume that   ($X_t$,$\mathcal{F}_t$,$\mathbf{P}_{x}$) 
is a regular diffusion on $S$, where $S$  is a interval of
$\mathbb{R}$.  
We assume that the differential operator $L$ of the diffusion 
is given by 
\begin{equation} \label{E:operator}
        L= \frac{1}{2} \sigma^2(x) \frac{d^2}{dx^2} + 
        b(x) \frac{d}{dx}
\end{equation}
where $\sigma^2 \colon S  \mapsto \mathbb{R}^+$, $b \colon S \mapsto 
\mathbb{R}$  satisfy Condition~\ref{C:existence}.  Moreover we assume that
$b/\sigma^2 \in C([c,d]) \cap C^2((c,d)$  where $c ,d \in S$ and $c <d$.
It follows that $s \colon [c,d] \to \mathbb{R}$ defined by 
\begin{equation}
s(x)=\int_{c}^x \exp \{ - \int_c^y \frac{2 b(z)}{\sigma^2 (z)} \, dz \} \, dy
\end{equation}
belongs to $C([c,d]) \cap C^2((c,d))$.  It is an elementary 
exercise to verify that $s$ satisfies the equation $Ls=0,$ with 
initial condition $s(c)=0$.
It follows by Theorem 13.16 volume II of Dynkin \cite{eD65} that 
\begin{equation}\label{E:calc_mom}
\{f(d)-f(c)\}\frac{s(x)}{s(d)} + f(c)=
\mathbf{E}_{x} f(X_{\tau_c \wedge \tau_d}))  
\end{equation}
  Let $\Lambda$ be an interval of $\mathbb{R}$. 
Assume that 
$(X_t,\mathcal{F}_t,\mathbf{P}_{x}^{\lambda})_{\lambda \in \Lambda}$
is a parametric set of diffusions on $\mathbb{R}$, with sample space
$(\Omega,\mathcal{F}_{\infty})$ and differential operators 
$(L_{\lambda})_{\lambda \in \Lambda}$,
where $L_{\lambda}$ is given by the formula
\begin{equation} \label{E:oneparametrizacion}
        L_{\lambda}= \frac{1}{2} \sigma^2(x,\lambda) \frac{d^2}{dx^2} + 
        b(x,\lambda) \frac{d}{dx} 
\end{equation}
Here we assume that $b(\cdot ,\lambda)$ and 
$\sigma^2(\cdot ,\lambda)$ satisfy the hypotheses
of this appendix where as before $c < d$ belongs to $S$ are fixed constants.
We wish to find conditions on $b, \sigma^2$ to guarantee that the
function given by $\lambda \mapsto  \mathbf{E}^{\lambda}_{x} 
f(X_{\tau_c \wedge \tau_d})$ is monotone decreasing (or increasing).  
For this we prove the following lemma:\\

\begin{lemma} \label{L:condmon}
Let $c<d$ be real numbers and let $\Lambda$ be a closed  interval of $\mathbb{R}$.  Let 
$f \colon [c,d] \times \Lambda \mapsto (0,\infty)$ be a jointly continuous
positive function such that $\partial f /\partial \lambda$ is also jointly 
continuous.  Let us assume that 
\begin{equation*}
\frac{\partial f}{\partial \lambda}(x,\lambda)=f(x,\lambda) g(x,\lambda)
\end{equation*}  
where  $g$ is a strictly increasing (strictly decreasing)
function in $x$ for each $\lambda \in \Lambda$.
It follows that  the function 
\begin{equation*}
h(x,\lambda)=\frac{\int_x^d f(y,\lambda) \, dy}{\int_c^x f(y,\lambda) \, dy}
\end{equation*}
is a strictly increasing (strictly decreasing) function in 
$\Lambda$ for any $x \in [c,d]$.
\begin{proof}
We prove $g$ strictly increasing in $x$ implies that $h$ 
is a strictly increasing function in $\lambda$.
(The proof of $g$ strictly decreasing implies that $h$ is 
strictly decreasing is similar.)  By the dominated convergence theorem
\small
\begin{equation*}
\frac{\partial h}{\partial \lambda}(x, \lambda)=
\frac{\int_c^x f(y,\lambda) \, dy
        \int_x^d \frac{\partial f}{\partial \lambda}(y,\lambda) \, dy -
        \int_x^d f(y,\lambda) \, dy  \int_c^x \frac{\partial f}
        {\partial \lambda}(y, \lambda)
        \, dy}{\{\int_c^x f(y, \lambda) \, dy\}^2}
\end{equation*}
\normalsize
The result follows  using the following inequalities:
\begin{align*}
(\int_x^d f(y, \lambda) \, dy ) g(x, \lambda) &<
        \int_x^d \frac{\partial f}{\partial \lambda}(y,\lambda) \, dy  &<  
                (\int_x^d f(y, \lambda) \, dy ) g(d, \lambda)\\
(\int_c^x f(y, \lambda) \, dy ) g(c, \lambda) &<
        \int_c^x \frac{\partial f}{\partial \lambda}(y, \lambda) \, dy &<
                (\int_c^x f(y, \lambda) \, dy ) g(x, \lambda)   
\end{align*}
\end{proof}
\end{lemma}

\begin{corollary} \label{C:criteria}
  Let $(X_t,\mathcal{F}_t,\mathbf{P}_{x}^{\lambda})_{\lambda \in \Lambda}$
be a parametric set of diffusions with
differential operators $(L_{\lambda})_{\lambda \in \Lambda}$ as in equation~(\ref{E:oneparametrizacion}). Assume that $b( \cdot,\lambda)$ 
and $\sigma^2(\cdot ,\lambda)$ 
satisfy the hypotheses
of this appendix where as before $c < d$ in  $S$ are fixed constants.
Let $\eta \colon [c,d] \times \Lambda \mapsto \mathbb{R}$ be defined by the
formula $\eta(x,\lambda)=\mathbf{E}^{\lambda}_{x} f(X_{\tau_c \wedge \tau_d})$.
If  $\partial/\partial \lambda (b/\sigma^2) (x,\lambda)>(<)0$ 
for all $(x,\lambda) \in [c,d] \times \Lambda$ then $\eta(x,\lambda)$ 
is a strictly decreasing (increasing) function on $\lambda$ for all 
$x \in [c,d]$
\end{corollary}
\begin{proof}
The result follows by Lemma~\ref{L:condmon} and equation~(\ref{E:oneparametrizacion}).
\end{proof}
Finally we mention a result that allows us, in the case of diffusions 
with an one-dimensional state space, to compute the expected values
of exit times from open sets.\\

\begin{lemma} \label{L:avwaiting}
Let  ($X_t$,$\mathcal{F}_t$,$\mathbf{P}_{x}$) be a regular diffusion on 
$S$, where $S$  is an interval of $\mathbb{R}$.  We assume that $L$ 
is the differential
operator of the diffusion where $L$ is defined
by equation~(\ref{E:operator}) and we assume that $\sigma^2$ and 
$b$ satisfy the hypothesis  of this appendix. Set
\begin{equation*} 
\varphi(x)=\exp \{- \int_c^x \frac{2b(z)}{\sigma^2(z)} \, dz\}
\end{equation*}
Then
\begin{multline} \label{E:avwaiting2}
\mathbf{E}_x \tau_c \wedge \tau_d =\eta(x)=
        -\int_c^x 2 \varphi(y) \int_c^y \frac{1}{\sigma^2(z)
        \varphi(z)}\, dz \,dy +\\
         \frac{\int_c^x \varphi(z) \, dz}{\int_c^d \varphi (z) dz}
        \int_c^d 2 \varphi(y) \int_c^y \frac{1}{\sigma^2(z) \varphi(z)} \, dz \, dy
\end{multline}
Moreover $u(x)=\mathbf{E}_x (\tau_c \wedge \tau_d)^2 < \infty$ 
for any $x \in [c,d]$ and $u$ is the solution of the differential equation
\begin{equation} \label{E:avwaiting3}
Lu=-\eta
\end{equation}
with boundary conditions $u(c)=u(d)=0$
\end{lemma}
\begin{proof}
Equation~(\ref{E:avwaiting2}) follows by Theorem 13.16 
volume II of Dynkin \cite{eD65} and a straightforward
computation.  The later part of the lemma follows by Theorem I.15.3 of
Gihman and A.V. Skorohod \cite{GS72}.
\end{proof}

\end{document}